\documentclass[12pt]{article}
\usepackage{epsfig}
\usepackage{amsmath}
\usepackage{amssymb}

\newtheorem{thm}{Theorem}[section] 
\newtheorem{lemma}[thm]{Lemma} 
\newtheorem{Def}[thm]{Definition} 
\newtheorem{cor}[thm]{Corollary} 
\newtheorem{prop}[thm]{Proposition} 
\newtheorem{rmk}[thm]{Remark} 

\newtheorem{examp}[thm]{Example}
\newtheorem{hyp}[thm]{Hypothesis}

\newcommand{\proof}{\noindent {\bf Proof} \hspace{0.2in}} 
\newcommand{\qed}{\hfill\mbox{\raggedright\rule{.07in}{.1in}}
  \vspace{1ex}} 
\newcommand{\dps}{\displaystyle}
\newcommand{\Section}[1]{\section{#1} \setcounter{equation}{0}}

\title{Toroidal normal forms for bifurcations in retarded functional
differential equations I: Multiple Hopf and transcritical/multiple
Hopf interaction}
 
\author{Younsun Choi\\Department of Mathematics and Statistics\\
University of Ottawa\\Ottawa, ON K1N 6N5\\CANADA  
\and Victor G. LeBlanc\\Department of
  Mathematics and Statistics\\University of Ottawa\\Ottawa, 
ON K1N 6N5\\CANADA}
\date{May 18, 2005}

\begin{document}

\maketitle

\begin{abstract}
For finite-dimensional bifurcation problems, it is
well-known that it is possible to compute normal forms which possess nice symmetry
properties.  Oftentimes, these symmetries may allow for a partial decoupling of
the normal form into a so-called ``radial'' part and an ``angular''
part.  Analysis of the radial part usually gives an enormous
amount of valuable
information about the bifurcation and its unfoldings.  In this paper,
we are interested in the case where such bifurcations occur in
retarded functional differential equations, and we revisit the
realizability and restrictions problem for {\em the class of radial equations} by
nonlinear delay-differential equations.  Our analysis allows us to
recover and considerably generalize recent results by Faria and
Magalh$\tilde{\mbox{\rm a}}$es \cite{FM96} and by Buono and B\'elair \cite{BB}.
\end{abstract}

\pagebreak
\Section{Introduction}
Delay-differential equations are used extensively in the modeling
of a multitude of phenomena in the life sciences \cite{BBL,Kuang,LM}, physics
\cite{HFEKGG,VTK}, atmospheric sciences \cite{SS}, engineering \cite{SC},
economics \cite{ZW} and beyond.  
This has motivated a flurry of activity
on the mathematical side to try to understand the behavior of this
class of equations and to
develop a theoretical framework suitable for their analysis.

It is now well-understood that retarded functional differential
equations (RFDEs), a class which contains delay-differential equations,
behave for the most part like infinite-dimensional ordinary
differential equations.  The upshot is that many of the techniques
and theoretical results of finite-dimensional geometrical dynamical systems are portable
to RFDEs.  In particular, versions of the stable/unstable and center
manifold theorems in neighborhoods of an equilibrium point exist for
RFDEs \cite{HL}.
For example, near a bifurcation point in a RFDE, the flow is essentially governed
by a vector field on an invariant center manifold.  This has allowed for
the successful application of the vast machinery of bifurcation theory
to many problems which are modeled by RFDEs, e.g. \cite{BC94,SC}.  
Parallel to this, techniques for simplifying vector fields via normal
form changes of coordinates have been adapted to RFDEs \cite{FMTB,FMH}, and has
allowed for further insight into the qualitative behavior of RFDEs.

This paper is concerned with the bifurcation theory
of RFDEs.  In particular, we will be interested in the so-called {\em
  realizability problem} for normal forms of vector fields which arise
via center manifold reduction of RFDEs.

\vspace*{0.15in}
\noindent
{\em Realizability problem:}

\vspace*{0.15in}
\noindent
Suppose $B$ is an arbitrary $m\times m$ matrix.  For the
sake of simplicity, suppose additionally that all eigenvalues of $B$
are simple.  
Let $C([-r,0],\mathbb{R})$ designate the space of continuous functions
from the interval $[-r,0]$ into $\mathbb{R}$, and for any continuous
function $z$, define $z_t\in
C([-r,0],\mathbb{R})$ as $z_t(\theta)=z(t+\theta)$, $-r\leq\theta\leq 0$.
It is then possible \cite{FMR} to construct a bounded linear operator 
$L:C([-r,0],\mathbb{R})\longrightarrow\mathbb{R}$ such that the
infinitesimal generator $A_0$ for the flow associated with the functional
differential equation 
\begin{equation}
\dot{z}(t)=L\,z_t
\label{linfde1}
\end{equation}
 has a spectrum which
contains the eigenvalues of $B$ as a subset. 
Thus, there exists an $m$-dimensional subspace $P$ of
$C([-r,0],\mathbb{R})$ which is invariant for the flow generated by
$A_0$, and the flow on $P$ is given by the linear ordinary
differential equation (ODE)
\[
\dot{x}=Bx.
\]
In our case, we will be especially interested in the case where the
eigenvalues of $B$ all have zero real parts, and the spectrum of
$A_0$ does not contain any elements with zero real part other than
those which belong to the spectrum
of $B$.

Now, suppose (\ref{linfde1}) is modified by the addition of a
nonlinear delayed term
\begin{equation}
\dot{z}(t)=L\,z_t+az(t+\tau)^2,
\label{nonlinfde1}
\end{equation}
where $a\in\mathbb{R}$ is some coefficient and $\tau\in [-r,0]$ is the
delay time.  Then the center manifold theorem for RFDEs\cite{HL} can be used to
show that the flow for (\ref{nonlinfde1}) admits an $m$-dimensional
locally invariant center manifold on which the dynamics associated
with (\ref{nonlinfde1}) are given by a vector field which, to
quadratic order, is of the form
\begin{equation}
\dot{x}=Bx+ag(x),
\label{realizeode1}
\end{equation}
where $g:\mathbb{R}^m\longrightarrow\mathbb{R}^m$ is a fixed homogeneous
quadratic polynomial which is completely determined by $L$ and $\tau$,
and $a$ is the same coefficient which appears in (\ref{nonlinfde1}).
We immediately notice that for fixed $L$ and $\tau$,
(\ref{realizeode1}) has at most one degree of freedom in the quadratic
term, corresponding to the one
degree of freedom in the quadratic term in (\ref{nonlinfde1}).
However, whereas one degree of freedom is sufficient to
describe the general scalar quadratic term involving one delay in
(\ref{nonlinfde1}), it is largely insufficient (if $m>1$) to describe
the general homogeneous quadratic polynomial
$f:\mathbb{R}^m\longrightarrow\mathbb{R}^m$.
Therefore, there exist $m$-dimensional vector fields $\dot{x}=Bx+f(x)$
(where $f$ is homogeneous quadratic) which can not be realized by
center manifold reduction (\ref{realizeode1}) of any RFDE of the form (\ref{nonlinfde1}).
One quickly notices that the situation could be improved if we allow the
nonlinear terms in (\ref{nonlinfde1}) to depend on more than one
delayed times, i.e.
\begin{equation}
\dot{z}(t)=L(z_t)+\sum_{\stackrel{i_1,\ldots,i_j=0}{i_1+\cdots+i_j=2}}^2
  a_{i_1i_2\cdots
  i_j}(z(t+\tau_1))^{i_1}\cdots (z(t+\tau_j))^{i_j},
\label{nonlinfde2}
\end{equation}
where the $a_{i_1i_2\cdots i_j}$ are real coefficients and
$\tau_1,\ldots,\tau_j\in [-r,0]$ are the delay times.
The center manifold equations for (\ref{nonlinfde2}) truncated to
quadratic order are
\begin{equation}
\dot{x}=Bx+\sum_{\stackrel{i_1,\ldots,i_j=0}{i_1+\cdots+i_j=2}}^2 a_{i_1i_2\cdots
  i_j}\,g_{i_1i_2\cdots i_j}(x),
\label{realizeode2}
\end{equation}
where $g_{i_1i_2\cdots i_j}:\mathbb{R}^m\longrightarrow\mathbb{R}^m$
are fixed homogeneous quadratic polynomials which are completely
determined by $L$ and $\tau_1,\ldots,\tau_j$.  Thus, the subspace of $m$-dimensional vector fields $\dot{x}=Bx+f(x)$
(where $f$ is a homogeneous quadratic) of the form (\ref{realizeode2})
is potentially larger
than those of the form (\ref{realizeode1}).  Of course, there is
nothing particularly special about the quadratic order, and one could
repeat the above discussion to include progressively higher order
nonlinearities.  Without loss of generality, we could also limit our
attention to only those $f$ which are in normal form with respect to
the matrix $B$.
The particular version of the {\em realizability problem} which will
interest us in this paper is the following:

\vspace*{0.15in}
{\em
Given:
\begin{itemize}
\item an $m\times m$ matrix $B$ whose spectrum consists solely of simple
  eigenvalues with zero real parts,
\item a bounded linear operator 
$L:C([-r,0],\mathbb{R})\longrightarrow\mathbb{R}$ such that the
infinitesimal generator $A_0$ for the flow associated with the functional
differential equation (\ref{linfde1})
has a spectrum which
contains the eigenvalues $B$ as a subset, and no other part of its
spectrum on the imaginary axis 
\item an integer $\ell\geq 2$
\item a polynomial $f:\mathbb{R}^m\longrightarrow\mathbb{R}^m$ of
  degree $\ell$ such
  that $f(0)=0$ and $Df(0)=0$, and $f$ is in normal form with respect
  to the matrix $B$
\end{itemize}
does there exist an RFDE of the form
\begin{equation}
\dot{z}(t)=L\,z_t+F(z(t+\tau_1),\ldots,z(t+\tau_j)),
\label{nonlinfde3}
\end{equation}
such that the center manifold equations for (\ref{nonlinfde3}), in
normal form and
truncated to order $\ell$, are $\dot{x}=Bx+f(x)$?
}

\vspace*{0.15in}
This question was answered in the affirmative for scalar RFDEs \cite{H85,H86} and in \cite{FMR} for
$n$-dimensional RFDEs in general.  In the scalar case (which will be
of interest to us), the result states that there is generically a solution to the realizability problem
for general $f$ as stated above if $j$ (the number of distinct delays
in (\ref{nonlinfde3})) is at least
equal to $m$ (the dimension of the center subspace).  

The main purpose of this paper is related to the optimality of the above
sufficient number ($j=m$) of delays, in light of some recent results
by Faria and Magalh$\tilde{\mbox{\rm a}}$es \cite{FM96}, and by Buono
and B\'elair \cite{BB}, which we now describe.

\vspace*{0.15in}
\noindent
{\em Simple Hopf, $(0,\pm\,i\omega)$, and $(\pm\,i\omega_1,\pm\,i\omega_2)$ bifurcations:}

\vspace*{0.15in}
\noindent
In \cite{FM96} and \cite{BB}, the authors consider the optimality of
the solution $j=m$ to the
realizability problem for scalar RFDEs in some special cases.
Consider one of the following three separate cases for the matrix $B$:
\begin{itemize} 
\item $B$ is a 
$2\times 2$ matrix whose eigenvalues are
$\pm\,i\omega$, $\omega>0$,
\item $B$ is a $3\times 3$ matrix whose eigenvalues are $0$ and
  $\pm\,i\omega$, $\omega>0$,
\item $B$ is a $4\times 4$ matrix whose eigenvalues are
  $\pm\,i\omega_1$, $\pm\,i\omega_2$, where $\omega_1>0$ and $\omega_2>0$
  are rationally incommensurate.
\end{itemize}
Let
$L:C([-r,0],\mathbb{R})\longrightarrow\mathbb{R}$ be a bounded linear
operator such that the
infinitesimal generator $A_0$ for the flow associated with the functional
differential equation (\ref{linfde1})
has a spectrum which
contains the eigenvalues of $B$ as a subset, and has no other part of
its spectrum on the imaginary axis.
Therefore, a general RFDE of the form
\begin{equation}
\dot{z}(t)=L\,z_t+N(z_t),
\label{nonlinfde4}
\end{equation}
where $N(0)=0$, $DN(0)=0$ has an equilibrium solution $z=0$
undergoing respectively a simple Hopf bifurcation, a
steady-state/Hopf interaction, or a non-resonant double Hopf bifurcation.  Normal forms and versal unfoldings for
each of these bifurcations are well-known.  In the first case, using
normal form changes of coordinates and then converting to
polar coordinates $x_1=\rho\cos\,\theta$, $x_2=\rho\sin\,\theta$ for the center manifold, the normal form (to cubic
order) is
\begin{equation}
\dot{\rho}=a\rho^3,\,\,\,\,\,\,\,\,\,\dot{\theta}=\omega+b\rho^2.
\label{hopfintro1}
\end{equation}
If the coefficient $a$ in (\ref{hopfintro1}) is non-zero, then the higher-order terms have
no qualitative effects.  In this case, the $\dot{\rho}$ equation
completely determines the bifurcation.  Now, from the above-mentioned
solution $j=2$ to
the realizability problem \cite{H85,H86}, we can conclude that if $N(z_t)$ is of the
form $N(z_t)=F(z(t+\tau_1),z(t+\tau_2))$, then any value of $a$ and $b$
in
(\ref{hopfintro1}) can be realized by means of center manifold
reduction of (\ref{nonlinfde4}).  However, Faria and
Magalh$\tilde{\mbox{\rm a}}$es
show that, in fact, any value of $a$ in the determining $\dot{\rho}$ equation of (\ref{hopfintro1}) can be
generically realized if $N(z_t)$ {\em only involves one delay},
i.e. $N(z_t)=F(z(t+\tau_1))$.  Similarly, they show that the versal
unfolding
\[
\dot{\rho}=\lambda\rho+a\rho^3
\]
of the $\dot{\rho}$ equation in (\ref{hopfintro1}) can be realized by
a RFDE of the form
$\dot{z}(t)=L(z_t)+\lambda\,z(t+\tau_1)+F(z(t+\tau_1))$.

For the $(0,\pm\,i\omega)$ case, the center manifold is
three-dimensional.  Normal form changes of coordinates
and the use of center manifold coordinates $x_1$, $x_2=\rho\cos\,\theta$,
$x_3=\rho\sin\,\theta$ yield the following equations (to quadratic order)
\[
\begin{array}{lll}
\dot{x}_1&=&b_{1}x_1^2+b_{2}\rho^2\\[0.15in]
\dot{\rho}&=&a_1x_1\rho\\[0.15in]
\dot{\theta}&=&\omega+O(|x_1,\rho|^2).
\end{array}
\]
If the coefficients $b_1$, $b_2$ and $a_1$ satisfy certain generic
non-degeneracy conditions, then the higher-order terms have no
qualitative effects.  The bifurcation is thus characterized by the
$\dot{x}_1$, $\dot{\rho}$ subsystem.  In this case, Faria and
Magalh$\tilde{\mbox{\rm a}}$es
show that any value of $b_1$, $b_2$ and $a_1$ can be realized 
by center manifold reduction of an RFDE involving only 2
delays (which is less than the predicted value ($j=m=3$) from the solution to the
realizability problem \cite{H85,H86}).

Finally, a similar result was shown for the non-resonant double Hopf
bifurcation in \cite{BB}: whereas the solution of the realizability
problem \cite{H85,H86}
predicts that $j=m=4$ delays are sufficient, it is shown in \cite{BB} that
2 delays are sufficient to realize, to cubic order, the ``radial
part'' ($\dot{\rho_1}$, $\dot{\rho_2}$) of the center manifold equations
\[
\begin{array}{lll}
\dot{\rho_1}&=&(\mu_1+a_{11}\rho_1^2+a_{12}\rho_2^2)\rho_1\\[0.15in]
\dot{\rho_2}&=&(\mu_2+a_{21}\rho_1^2+a_{22}\rho_2^2)\rho_2\\[0.15in]
\dot{\theta_1}&=&\omega_1+O(|\mu_1,\mu_2|,|\rho_1,\rho_2|^2)\\[0.15in]
\dot{\theta_2}&=&\omega_2+O(|\mu_1,\mu_2|,|\rho_1,\rho_2|^2).
\end{array}
\]

In all three cases above, it is also shown that this smaller
number of delays (1 in the simple Hopf case and 2 in both the
$(0,\pm\,i\omega)$ and $(\pm\,i\omega_1,\pm\,i\omega_2)$ cases) is
optimal, in the sense that anything less will lead to restrictions on
realizability of the various coefficients which appear in these normal
forms, and consequently to restrictions on
the possible phase portraits in the classification of the versal
unfolding of these respective singularities.

\vspace*{0.15in}
\noindent
{\em Overview:}

\vspace*{0.15in}
\noindent
The questions of realizability and restriction for normal forms and
unfoldings of bifurcations in RFDEs are particularly important from a modeling
point of view.  Indeed, given a specific RFDE model (perhaps depending
on many parameters) undergoing a local bifurcation, it is important
to be able to characterize the range of possible dynamics accessible
from within the model near the bifurcation point.  From our discussion above, we see that knowledge of
the abstract finite-dimensional bifurcation problem and its unfoldings is not
sufficient in general to answer this question.  
Indeed, the specific form of the RFDE may restrict this range of
possible dynamics.
Depending on the functional form
of the RFDE (e.g. how many distinct delays are involved), some phase
diagrams which are possible in the unfolding of this given bifurcation may not be realizable in
the RFDE.  This could have important consequences in the
interpretation of the model, especially as it pertains to the actual
phenomenon being modeled.

The purpose of this paper is to further study these issues of
realizability and restrictions, in light of the
previously discussed results in
\cite{FM96} and \cite{BB}.  We will develop a
unified theoretical framework for these results, which will
consequently allow for
considerable generalizations of these
results.

Specifically, we will exploit and generalize the following common elements of
the three specific cases studied in \cite{FM96} and \cite{BB}: it is possible to
make a canonical choice of normal form transformations on the center
manifold equations which lead to a normal form with nice symmetry
properties -- it is equivariant with respect to an action of a torus
group.  This toroidal equivariance can be used to achieve a partial
decoupling of the normal form into a ``radial part'' and an ``angular
part''.  In many cases (in particular, in the three specific cases
studied in \cite{FM96} and \cite{BB}), the radial part characterizes the essential
features of the bifurcation.  It is then reasonable to

\vspace*{0.1in}
\noindent
{\em investigate the realizability problem for the radial part of the
normal form and its unfoldings}, 

\vspace*{0.1in}
\noindent
which is the goal we seek.  Along with developing a
theoretical framework to achieve this goal, we will in fact make the
following generalizations to
the results of \cite{FM96} and \cite{BB}:
\begin{itemize}
\item we will assume that the spectrum of the matrix $B$ consists of
  simple eigenvalues, and has one of the
  two following forms
\begin{equation}
\mbox{\rm
  spec}(B)=\{\pm\,i\omega_1,\ldots,\pm\,i\omega_p\}\,\,\,\,\,\,\mbox{\rm
  or}\,\,\,\,\,\,
\mbox{\rm spec}(B)=\{0,\pm\,i\omega_1,\ldots,\pm\,i\omega_p\},
\label{specBhypintro}
\end{equation}
where $\omega_1,\ldots,\omega_p>0$ are
independent over the rationals,
\item we will investigate realizability of the radial part of the
  normal form {\em to any order}, and not just quadratic or cubic.
  This is important in cases where nonlinear degeneracies are present.
\end{itemize}

Similarly to \cite{BB,FM96},
we will limit our analysis to the case of {\em scalar} RFDEs.
While studying the realizability problem in the context of
general $n>1$ dimensional systems of RFDEs is certainly very
important, our computations indicate that there is an enormous
increase in algebraic complexity involved.  Consequently, a unified concise
framework allowing for the simultaneous analysis of all cases (scalar and systems) 
appears at this point to be a difficult, albeit not impossible, goal
to achieve.  

In the second case for the spectrum of the matrix $B$ in
(\ref{specBhypintro}), there is an technical subtlety which
arises as it pertains to the possible unfoldings of this singularity.
In fact, there are two algebraically different ways to construct an
unfolding, depending on whether the ``steady-state'' mode
(correspongind to the 0 eigenvalue) in the
interaction is of saddle-node type or of transcritical type.  It turns
out that the transcritical case can be treated in the same framework as the
first case for $\mbox{\rm spec}(B)$ in (\ref{specBhypintro}), but it would
be extremely cumbersome to attempt to treat the
saddle-node case within this same framework.  Therefore, we have
chosen to treat the saddle-node case in a separate paper \cite{CL}.

This paper is organized as follows.  In Section 2, we will give a
brief review of the theory of the center manifold reduction and normal
form transformations of RFDEs as developed by Faria and 
Magalh$\tilde{\mbox{\rm a}}$es in \cite{FMTB,FMH}.  In Section 3, we will review
how to make a canonical choice of normal form which possesses useful toroidal
symmetry properties, and exploit this symmetry to achieve a partial
decoupling of the normal form into a ``radial part'' and an ``angular
part''.  The radial part possesses residual reflectional-type
symmetry, which will be important for the subsequent analysis.  In
Section 4, we set a framework and establish an important surjectivity
result which will be crucial to study the realizability
problem for the radial part.  Our main results on realizability and
restrictions are given in Section 5.  In this section, we also give
some results which hint at the rudiments of a singularity theory for RFDEs.
In Section 6 we show how the specific results of \cite{FM96} and \cite{BB} are recovered by our
main results.  We end with some concluding remarks in Section 7.  Some
of the proofs in Section 3 are relegated to the appendices.

\Section{Functional Analytic Framework}
In this section we will briefly recall some standard results and
terminology in the bifurcation theory of RFDEs in order
to establish the notation.  For more details, see \cite{FMTB,FMH,HL}.

\subsection{Infinite dimensional parameterized ODE}
Suppose $r>0$ is a given real
number and $C=C\left(  \left[  -r,0\right]  ,\mathbb{R}\right)  $ is the
Banach space of continuous functions from $\left[  -r,0\right]  $ into
$\mathbb{R}$ with supremum norm.  We define $z_{t}\in C$ as $z_{t}\left(  \theta\right)
=z\left(  t+\theta\right)  ,-r\leq\theta\leq 0.$ 
Let us consider the
following parameterized family of nonlinear retarded functional differential equations
\begin{equation}
\dot{z}\left(  t\right)  =L(\mu)z_{t}+F\left(  z_{t},\mu\right)  ,
\label{y1}%
\end{equation}
\noindent where $L:C\times\mathbb{R}^{s}\rightarrow\mathbb{R}$ is
a parameterized family of bounded linear
operators from $C$ into $\mathbb{R}$ and $F$ is
a smooth function from $C\times\mathbb{R}^{s}$ into $\mathbb{R}$.
In this paper, we will assume the following hypothesis on $F$
\begin{hyp}
$F(0,0)=0$ and
$DF(0,0)=0$.
\label{Fhyp}
\end{hyp}
A consequence of Hypothesis \ref{Fhyp} is that in a Taylor expansion
of (\ref{y1}), there
are no terms which are $z_t$ independent and linear in $\mu$.  
While this is not a restriction in one of the cases we
will be studying in this paper (multiple non-resonant Hopf
bifurcation), it is a restriction in the other case
(steady-state/multiple non-resonant Hopf interaction).  
Note however that Hypothesis \ref{Fhyp} includes as a special case the physically
interesting case in which $z=0$ is an equilibrium for (\ref{y1}) for
all $\mu$, which is the case of interaction between a
transcritical bifurcation and multiple non-resonant Hopf bifurcation
(see \cite{L79}).
In a sequel to this paper, we will relax Hypothesis \ref{Fhyp} to simply
$F(0,0)=0$ and $D_1F(0,0)=0$, which is the generic saddle-node case.  This relaxation of
Hypothesis \ref{Fhyp} leads to some technical complications which would make a
unified treatment of both cases simultaneously extremely cumbersome
and lengthy.
Therefore, for the sake of clarity, we have decided to treat these two
cases separately (see \cite{CL}).
  
The bounded linear map $L(\mu)$
can be represented in an integral form as
\[
L(\mu)\phi=\int_{-r}^{0}\left[  d\eta_{\mu}\left(\theta\right)  \right]  \phi\left(
\theta\right)  ,
\]
\noindent where $\eta_{\mu}\left(\theta\right)  $ is a measurable
function on $\left[  -r,0\right]  $.  Denote $L_0\equiv
L(0)$, and rewrite (\ref{y1}) as
\begin{equation}
\dot{z}\left( t\right) = 
L_0z_t+(L(\mu)-L_0)z_t+F(z_t,\mu)=L_0z_t+\widetilde{F}(z_t,\mu),
\label{y1p}
\end{equation}
where $\widetilde{F}(z_t,\mu)=(L(\mu)-L_0)z_t+F(z_t,\mu)$.

Let $A(\mu)$ be the infinitesimal generator of the flow for the linear
system $\dot{z}=L(\mu)z_t$, with spectrum $\sigma(A(\mu))$, and denote by
$\Lambda_{\mu}$ the set of eigenvalues of $\sigma(A(\mu))$ with zero
real part.

The set $\Lambda_0$, which consists of the roots of the characteristic equation
\begin{equation}
\det\Delta\left(  z\right)  =0,\text{ \ \ \ }\Delta\left(  z\right)
=z-\int\limits_{-r}^{0}\left[  d\eta_0\left(\theta\right)  \right]
e^{z\theta}, \label{ychar}%
\end{equation}
with zero real part
will play an important role.
\begin{hyp}
Throughout the rest of the paper, we assume the following hypotheses on $\Lambda_{\mu}$ and $\Lambda_0$:
\begin{enumerate}
\item[(a)]
$\mbox{\rm Card}(\Lambda_{\mu})<\mbox{\rm
    Card}(\Lambda_0)\,\,\,\,\mbox{\rm for $\mu$ small}$,
\item[(b)] Each element of $\Lambda_0$ is a simple eigenvalue of
  $A(0)$, and $\Lambda_0$ has one of the following two forms:
\[
\begin{array}{l}
\Lambda_0=\{\pm\,i\omega_1,\ldots,\pm\,i\omega_p\}\,\,\,\,\,\mbox{\rm
  (multiple non-resonant Hopf bifurcation), or}\\ \\
\Lambda_0=\{0,\pm\,i\omega_1,\ldots,\pm\,i\omega_p\}\,\,\,\,\,\mbox{\rm
  (steady-state/multiple non-resonant Hopf interaction),}
\end{array}
\]
where
$\omega_1,\ldots,\omega_p$, are independent over the rationals,
i.e. if $r_1,\ldots,r_p$ are rational numbers such that
$
\sum_{j=1}^p\,r_j\omega_j=0
$, then $r_1=\cdots=r_p=0$.
\end{enumerate}
\label{spectralhyp}
\end{hyp}

Let $P$ be the invariant subspace for $A_{0}\equiv A(0)$ associated with the
eigenvalues in $\Lambda_0$, and let $\Phi=\left(  \varphi_{1}\,\,\ldots\,\,
\varphi_{m}\right)$ be a matrix whose columns form a basis for $P$.  

In a similar manner, we can define an
invariant space, $P^{\ast},$ to be the generalized eigenspace of the
transposed system, $A_{0}^{T}$\ associated with $\Lambda_0,$ having as
basis the rows of the matrix
$\Psi=$\textrm{col}$\left(  \psi_{1},\ldots,\psi_{m}\right)$. Note that the
transposed system, $A_{0}^{T}$ is defined over a dual space $C^{\ast}=C\left(
\left[  0,r\right]  ,\mathbb{R}\right),$ and each element of $\Psi$
is included in $C^{\ast}.$ The bilinear form between $C^{\ast}$ and $C$ is
defined as
\begin{equation}
\left(  \psi,\phi\right)  =\psi\left(  0\right)  \phi\left(  0\right)
-\int\limits_{-r}^{0}\int\limits_{0}^{\theta}\psi\left(  \zeta-\theta\right)
\text{ }\left[  d\eta_{0}\left(  \theta\right)  \right]  \text{ }\phi\left(
\zeta\right)  \text{ }d\zeta. \label{y11}%
\end{equation}
Note that $\Phi$ and $\Psi$ satisfy $\dot{\Phi}=B\Phi,$
$\dot{\Psi}=-\Psi B,$ where $B$ is an $m\times m$ matrix
whose spectrum coincides with $\Lambda_0$.

We can normalize
$\Psi$ such that $\left(  \Psi
,\Phi\right)  =I$, and we can decompose the space $C$ using the splitting
$C=P\oplus Q$, where the complimentary space $Q$ is also invariant for $A_0$.

Faria and Magalh$\tilde{\mbox{\rm a}}$es \cite{FMTB,FMH} show that (\ref{y1}) can be written
as an infinite dimensional ordinary differential equation on the
Banach space 
$BC$ of functions from $[-r,0]$ into ${\mathbb R}$ 
which are uniformly continuous on $[-r,0)$ and with a jump discontinuity 
at $0$, using a procedure that we will now outline.
Define $X_0$ to be the function
\[
X_0(\theta)=\left\{\begin{array}{lc}
1&\theta=0\\[0.15in]
0&-r\leq\theta<0,
\end{array}\right.
\]
then the elements of $BC$ can be written as $\xi=\varphi+X_0\lambda$, with
$\varphi\in C$ and $\lambda\in {\mathbb R}$, so that
$BC$ is identified with $C\times {\mathbb R}$.  

Let $\pi:BC\longrightarrow P$ denote the projection
\[
\pi(\varphi+X_0\lambda)=\Phi [(\Psi,\varphi)+\Psi(0)\lambda],
\]
where $\varphi\in C$ and $\lambda\in {\mathbb R}$.  
We now decompose $z_t$ in (\ref{y1}) according to the splitting
\[
BC=P\oplus\mbox{\rm ker}\,\pi,
\]
with the property that $Q\subsetneq\,\mbox{\rm ker}\,\pi$, 
and get the following infinite-dimensional ODE system which is equivalent to (\ref{y1}):
\begin{equation}
\begin{array}{rcl}
\dot{x}&=&Bx+\Psi(0)\left[(L(\mu)-L_0)(\Phi\,x+y)+F(\Phi\,x+y,\mu)\right]\\[0.15in]
{\displaystyle\frac{d}{dt}\,y}&=&A_{Q^1}y+(I-\pi)X_0\,\left[(L(\mu)-L_0)(\Phi\,x+y)+F(\Phi\,x+y,\mu)\right],
\end{array}
\label{projfdep}
\end{equation}
where $x\in {\mathbb R}^m$,
$y\in Q^1\equiv Q\cap C^1$,
($C^1$ is the subset of $C$ consisting
of continuously differentiable functions),
and
$A_{Q^1}$ is the operator from 
$Q^1$ into 
$\mbox{\rm ker}\,\pi$ defined by
\[
A_{Q^1}\varphi=\dot{\varphi}+X_0\,[L\,\varphi-\dot{\varphi}(0)].
\]

\subsection{Faria and Magalh$\tilde{\mbox{\bf a}}$es normal form}

Consider the formal Taylor expansion of the nonlinear terms $\widetilde{F}$ in (\ref{y1p})
\[
\widetilde{F}(u,\mu)=\sum_{j\geq 2}\,\widetilde{F}_j(u,\mu),\,\,\,\,\,u\in\,C,\,\,\mu\in\mathbb{R}^s,
\]
where $\widetilde{F}_j(w)=H_j(w,\ldots,w)$, with $H_j$ belonging to the space of
continuous multilinear symmetric maps from
$(C\times\mathbb{R}^s)\times\cdots\times (C\times\mathbb{R}^s)$
($j$ times) to $\mathbb{R}$. 
If we denote $f_j=(f_j^1,f_j^2)$, where
\[
\begin{array}{rcl}
f_{j}^1(x,y,\mu)&=&\Psi(0)\,\widetilde{F}_j(\Phi\,x+y,\mu)\\[0.15in]
f_j^2(x,y,\mu)&=&(I-\pi)\,X_0\,\widetilde{F}_j(\Phi\,x+y,\mu),
\end{array}
\]
then (\ref{projfdep}) can be written as
\begin{equation}
\begin{array}{rcl}
\dot{x}&=&{\displaystyle Bx+\sum_{j\geq
    2}\,f_j^1(x,y,\mu)}\\[0.15in]
{\displaystyle\frac{d}{dt}\,y}&=&{\displaystyle A_{Q^1}y+\sum_{j\geq
    2}\,f_j^2(x,y,\mu)}
\end{array}
\label{y4}
\end{equation}

The spectral hypotheses we have specified in Hypothesis
\ref{spectralhyp} are sufficient to conclude that the non-resonance
condition of Faria and Magalh$\tilde{\mbox{\rm a}}$es \cite{FMTB,FMH} holds.  Consequently, using
successively at each order $j$ a near identity change of variables of the form
\begin{equation}
(x,y)=(\hat{x},\hat{y})+U_j(\hat{x},\mu)\equiv (\hat{x},\hat{y})+
(U^1_j(\hat{x},\mu),U^2_j(\hat{x},\mu)),
\label{nfcv}
\end{equation}
(where $U^{1,2}_j$ are homogeneous degree $j$ polynomials in
the indicated variables, with coefficients respectively in
$\mathbb{R}^m$ and $Q^1$)
system (\ref{y4}) can be put into formal normal form
\begin{equation}
\begin{array}{rcl}
\dot{x}&=&{\displaystyle Bx+\sum_{j\geq
    2}\,g_j^1(x,y,\mu)}\\[0.15in]
{\displaystyle\frac{d}{dt}\,y}&=&{\displaystyle A_{Q^1}y+\sum_{j\geq
    2}\,g_j^2(x,y,\mu)}
\end{array}
\label{y5}
\end{equation}
such that the center manifold is locally given by $y=0$ and the local
flow of (\ref{y1}) on this center manifold is given by
\begin{equation}
\dot{x}=Bx+\sum_{j\geq 2}\,g_j^1(x,0,\mu).
\label{y6}
\end{equation}
The nonlinear terms in (\ref{y6}) are in normal form in the
classical sense with respect to the matrix $B$.

\Section{Bifurcations with Toroidal Normal Forms}
With equation (\ref{y6}) in mind, in this section we will discuss
normal form transformations of the general parameterized system
\begin{equation}
\begin{array}{lll}
\dot{x}&=&Bx+f(x,\mu)\\
\dot{\mu}&=&0,
\end{array}
\label{y7}
\end{equation}
where the spectrum $\Lambda_0$
of the matrix $B$ is as in Hypothesis \ref{spectralhyp}(b).
As much as possible, we will treat both cases of Hypothesis \ref{spectralhyp}(b)
(i.e. whether or
not $\Lambda_0$ includes $0$) 
simultaneously by adopting a notation which uses integers $\kappa$ and
$d$, which should be interpreted as having the values $\kappa=2p$ and $d=p$
in the case where $\Lambda_0=\{\pm i\omega_1,\ldots,\pm i\omega_p\}$,
and
the values $\kappa=2p+1$ and $d=p+1$ in the case where
$\Lambda_0=\{0,\pm i\omega_1,\ldots,\pm i\omega_p\}$.

It will be
extremely useful to use complex coordinates for the last $2p$
components of the space $\mathbb{R}^{\kappa}$,
so that we can identify
\[
\mathbb{R}^{\kappa}=\left\{\begin{array}{lcl}
\{\,(x_1,\overline{x_1},\ldots,x_p,\overline{x_p})\,\,|\,\,x_j\in\mathbb{C},\,j=1,\ldots,p\,\}\,\,\,&\mbox{\rm
  if}&\,\,\kappa=2p\\[0.15in]
\{\,
(x_0,x_1,\overline{x_1},\ldots,x_p,\overline{x_p})\,\,|\,\,x_0\in\mathbb{R},\,x_j\in\mathbb{C},\,j=1,\ldots,p\,\}&\mbox{\rm
  if}&\kappa=2p+1.
\end{array}\right.
\]
Then, without loss of generality, we may assume that
\begin{equation}
B=\mbox{\rm
  diag}(i\omega_1,-i\omega_1,\ldots,i\omega_p,-i\omega_p)\,\,\,\,\,\mbox{\rm or}\,\,\,\,\,
B=\mbox{\rm
  diag}(0,i\omega_1,\-i\omega_1,\ldots,i\omega_p,-i\omega_p)
\label{BdiagHopf}
\end{equation}
depending on which case of Hypothesis \ref{spectralhyp}(b) is being considered. 

At times, it will be convenient to write (\ref{y7}) as
\begin{equation}
\dot{\tilde{x}}=\tilde{B}\tilde{x}+\tilde{f}(\tilde{x}),
\label{y7ext}
\end{equation}
where $\tilde{x}=(x,\mu)$, $\tilde{f}=(f,0)$ and
\begin{equation}
\tilde{B}=\left(\begin{array}{c|c}
B&0\\\hline
0&0\end{array}\right).
\label{tildeBform}
\end{equation}

The section is divided into four subsections.  
In the first, we will
give a brief review of results on symmetric normal forms with
parameters.  
Most (if not all) of these results are
largely well-known in the unparameterized case (see for example
\cite{ETCBI,GSSII}), 
and only minor modifications are required to obtain the parameterized
versions we present herein.

In the second subsection, we will define an {\em equivariant
  projection operator} which will be useful in the computation of
  symmetric normal forms.

In the third
subsection, we will specify how the symmetry of these normal forms can
be exploited in order to achieve a partial decoupling of the normal
form.

Finally, in the fourth subsection, we will introduce a splitting of
our spaces of polynomials which naturally decomposes any vector field into
a singular parameter independent part plus a perturbation.

\subsection{Normal forms and toroidal symmetry}
For $\tilde{B}$ as in (\ref{tildeBform}), let $\tilde{B}^t$ denote the
transpose of $\tilde{B}$ and
let $\Gamma=\overline{\{e^{s\tilde{B}^t}\,|\,s\in\mathbb{R}\}}$ (where the
closure is taken in the space of $(\kappa+s)\times (\kappa+s)$ matrices), and note
that $\Gamma$ is an abelian connected Lie group isomorphic
to $\mathbb{T}^p$, where
$\mathbb{T}^p$ is the
$p$-torus:
\begin{equation}
\mathbb{T}^p=\left\{\begin{array}{lcl}
\{\,\mbox{\rm
  diag}(e^{i\theta_1},e^{-i\theta_1},\ldots,e^{i\theta_p},e^{-i\theta_p},1,\ldots,1)\,\,\,|\,\,\,\theta_j\in\mathbb{S}^1,\,j=1,\ldots,p\,\}\,\,\,&\mbox{\rm if}&\,\,\kappa=2p\\[0.15in]
\{\,\mbox{\rm
  diag}(1,e^{i\theta_1},e^{-i\theta_1},\ldots,e^{i\theta_p},e^{-i\theta_p},1\ldots,1)\,\,\,|\,\,\,\theta_j\in\mathbb{S}^1,\,j=1,\ldots,p\,\}&\mbox{\rm if}&\kappa=2p+1.\end{array}\right.
\label{Tpdef}
\end{equation}
\begin{Def}
For a given integer $\ell\geq 2$, a given normed space $X$, and for
$\kappa=2p$ (respectively $\kappa=2p+1$), we denote
by $H^{\kappa+s}_{\ell}(X)$ the linear space of homogeneous polynomials
of degree $\ell$ in the $\kappa+s$ variables
$x=(x_1,\overline{x_1},\ldots,x_p,\overline{x_p})$ (respectively
$x=(x_0,x_1,\overline{x_1},\ldots,x_p,\overline{x_p})$) and
$\mu=(\mu_1,\ldots,\mu_s)$ with coefficients in $X$.
For $X=\mathbb{R}^{\kappa+s}$, define 
$H^{\kappa+s}_{\ell}(\mathbb{R}^{\kappa+s},\Gamma)\subset H^{\kappa+s}_{\ell}(\mathbb{R}^{\kappa+s})$
to be the
subspace of $\Gamma$-equivariant maps, i.e.
\[
\begin{array}{l}
\tilde{f}\in H^{\kappa+s}_{\ell}(\mathbb{R}^{\kappa+s},\Gamma)
\Longleftrightarrow\\[0.15in]
\tilde{f}\in H^{\kappa+s}_{\ell}(\mathbb{R}^{\kappa+s})\,\,\,\mbox{\rm and}\,\,\,
\gamma
\tilde{f}(\gamma^{-1}\tilde{x})=\tilde{f}(\tilde{x}),\,\,\forall\,\tilde{x}=(x,\mu)\in\mathbb{R}^{\kappa+s},\,\,\forall\,\gamma\in\Gamma.
\end{array}
\]
\end{Def}
For the general class of near-identity changes of variables
$\tilde{x}\mapsto \hat{x}+h(\hat{x})$ for (\ref{y7ext}), it is
well-known that we can eliminate from (\ref{y7ext}) all nonlinear
terms which are in the range of the
{\em homological operator}
\begin{equation}
\begin{array}{c}
{\cal L}_{\tilde{B}}:H^{\kappa+s}_{\ell}(\mathbb{R}^{\kappa+s})\longrightarrow H^{\kappa+s}_{\ell}(\mathbb{R}^{\kappa+s})\\[0.15in]
\tilde{f}\longmapsto ({\cal L}_{\tilde{B}}\tilde{f})(\tilde{x})=D\tilde{f}(\tilde{x})\tilde{B}\tilde{x}-\tilde{B}\tilde{f}(\tilde{x}).
\end{array}
\label{homdefine}
\end{equation}
Thus, we must define in
$H^{\kappa+s}_{\ell}(\mathbb{R}^{\kappa+s})$ a complimentary space to
$\mbox{\rm range}\,{\cal L}_{\tilde{B}}$.  Of course, such a space is
not unique.  However, there exists a nice
canonical choice which will be extremely useful for our purposes (see
for example \cite{ETCBI,GSSII}).
\begin{prop}
\[
H^{\kappa+s}_{\ell}(\mathbb{R}^{\kappa+s})=H^{\kappa+s}_{\ell}(\mathbb{R}^{\kappa+s},\Gamma)\oplus\mbox{\rm
  range}\,{\cal L}_{\tilde{B}}
\]
\label{prop_enf1}
\end{prop}
The usefulness of Proposition \ref{prop_enf1} is that it is
straightforward to compute the general element of $H^{\kappa+s}_{\ell}(\mathbb{R}^{\kappa+s},\Gamma)$.
\begin{lemma}
Let $\tilde{B}$ be as in (\ref{y7ext}).  Then a smooth vector field
$\tilde{f}:\mathbb{R}^{\kappa+s}\longrightarrow\mathbb{R}^{\kappa+s}$
is $\mathbb{T}^p$-equivariant if and only if
$\tilde{f}$ has one of the following forms
\begin{equation}
\tilde{f}(x,\mu)=
\left(
\begin{array}{c}
a_1(x_1\overline{x_1},\ldots,x_p\overline{x_p},\mu)\,x_1\\ \\
\overline{a_1(x_1\overline{x_1},\ldots,x_p\overline{x_p},\mu)\,x_1}\\ \\
\vdots\\ \\
a_p(x_1\overline{x_1},\ldots,x_p\overline{x_p},\mu)\,x_p\\ \\
\overline{a_p(x_1\overline{x_1},\ldots,x_p\overline{x_p},\mu)\,x_p},\\ \\
b_1(x_1\overline{x_1},\ldots,x_p\overline{x_p},\mu)\\ \\
\vdots\\ \\
b_s(x_1\overline{x_1},\ldots,x_p\overline{x_p},\mu),
\end{array}
\right)\,\,\,\mbox{\rm or}\,\,\,
\left(
\begin{array}{c}
a_0(x_0,x_1\overline{x_1},\ldots,x_p\overline{x_p},\mu)\\ \\
a_1(x_0,x_1\overline{x_1},\ldots,x_p\overline{x_p},\mu)\,x_1\\ \\
\overline{a_1(x_0,x_1\overline{x_1},\ldots,x_p\overline{x_p},\mu)\,x_1}\\ \\
\vdots\\ \\
a_p(x_0,x_1\overline{x_1},\ldots,x_p\overline{x_p},\mu)\,x_p\\ \\
\overline{a_p(x_0,x_1\overline{x_1},\ldots,x_p\overline{x_p},\mu)\,x_p},\\
\\
b_1(x_0,x_1\overline{x_1},\ldots,x_p\overline{x_p},\mu)\\ \\
\vdots\\ \\
b_s(x_0,x_1\overline{x_1},\ldots,x_p\overline{x_p},\mu),
\end{array}
\right)
\label{torus_nf}
\end{equation}
respectively if $\kappa=2p$ or $\kappa=2p+1$,
where $a_1,\ldots,a_p$ are smooth and complex-valued, and $a_{0},
b_1,\ldots b_s$ are
smooth and real-valued.  
\label{lemenf1}
\end{lemma}
\proof
This is a standard result which is a consequence
of Schwarz lemma \cite{Sch}.  See also \cite{GSSII}.
\hfill\qed

Proposition \ref{prop_enf1} and Lemma \ref{lemenf1} are not exactly in
a form suitable for our purposes, since the vector field $\tilde{f}$
in (\ref{y7ext}) has the special form $\tilde{f}=(f,0)$ which we require
our normal form changes of variables to preserve.  Since we are only interested in the first $\kappa$ components of
(\ref{y7ext}), we would like to obtain a splitting of
$H^{\kappa+s}_{\ell}(\mathbb{R}^{\kappa})$ akin to the splitting 
of $H^{\kappa+s}_{\ell}(\mathbb{R}^{\kappa+s})$
in Proposition \ref{prop_enf1}.
For this purpose, we will need the following
\begin{Def}
\hspace*{1in}
\newline
\begin{enumerate}
\item[(a)]
We define $H^{\kappa+s}_{\ell}(\mathbb{R}^{\kappa},\mathbb{T}^p)$ to be the
subset of $H^{\kappa+s}_{\ell}(\mathbb{R}^{\kappa})$ consisting of
mappings \mbox{\rm 
$f:\mathbb{R}^{\kappa+s}\longrightarrow\mathbb{R}^{\kappa}$}
whose components are of the form of the first $\kappa$ components of
(\ref{torus_nf}).
Note that $H^{\kappa+s}_{\ell}(\mathbb{R}^{\kappa},\mathbb{T}^p)$
consists precisely of the $\mathbb{T}^p$-equivariant elements of
$H^{\kappa+s}_{\ell}(\mathbb{R}^{\kappa})$; that is, 
\[
\begin{array}{l}
f\in
H^{\kappa+s}_{\ell}(\mathbb{R}^{\kappa},\mathbb{T}^p)\Longleftrightarrow\\
f\in H^{\kappa+s}_{\ell}(\mathbb{R}^{\kappa})\,\,\,\mbox{\rm
  and}\,\,\,
f(\gamma_0\,x,\mu)=\gamma_0\,f(x,\mu),\,\,\,\forall\gamma_0\in\Gamma_0,\,\,\forall\,(x,\mu)\in\mathbb{R}^{\kappa+s},
\end{array}
\]
where $\Gamma_0$ is the group of $\kappa\times\kappa$ matrices which
is isomorphic to $\mathbb{T}^p$, and is parameterized as
\begin{equation}
\Gamma_0=\left\{\begin{array}{lcl}
\{\,\mbox{\rm
  diag}(e^{i\theta_1},e^{-i\theta_1},\ldots,e^{i\theta_p},e^{-i\theta_p})\,\,\,|\,\,\,\theta_j\in\mathbb{S}^1,\,j=1,\ldots,p\,\}\,\,\,&\mbox{\rm if}&\,\,\kappa=2p\\[0.15in]
\{\,\mbox{\rm
  diag}(1,e^{i\theta_1},e^{-i\theta_1},\ldots,e^{i\theta_p},e^{-i\theta_p})\,\,\,|\,\,\,\theta_j\in\mathbb{S}^1,\,j=1,\ldots,p\,\}&\mbox{\rm if}&\kappa=2p+1,\end{array}\right.
\label{Tpdef2}
\end{equation}
\item[(b)] 
We define the following operator
\begin{equation}
\begin{array}{c}
{\cal L}_{B} : H^{\kappa+s}_{\ell}(\mathbb{R}^{\kappa})
\longrightarrow H^{\kappa+s}_{\ell}(\mathbb{R}^{\kappa})\\
f\longmapsto ({\cal
  L}_{B})(f)(x,\mu)=D_xf(x,\mu)Bx-Bf(x,\mu).
\end{array}
\label{homdefine2}
\end{equation}
Note that
${\cal L}_{\tilde{B}}(f,0)=({\cal L}_{B}f,0)$.
\end{enumerate}
\end{Def}
\begin{prop}
\[
H^{\kappa+s}_{\ell}(\mathbb{R}^{\kappa})=H^{\kappa+s}_{\ell}(\mathbb{R}^{\kappa},\mathbb{T}^p)\oplus\mbox{\rm
  range}\,
{\cal L}_{B}.
\]
\label{prop_enf2}
\end{prop}
\proof
The proof is given in the appendix.

\subsection{Equivariant projection}
In this section, we will construct an appropriate linear projection
associated with the splitting of
$H^{\kappa+s}_{\ell}(\mathbb{R}^{\kappa})$ given in Proposition
\ref{prop_enf2}.  This projection has very nice algebraic
properties, and will be useful when we prove our main results later.
\begin{Def}
Let ${\displaystyle\int_{\Gamma_0}\,d\gamma}$ denote the normalized
Haar integral on $\Gamma_0\cong\mathbb{T}^p$ (see (\ref{Tpdef2})).  We define the linear
operator
\[
\begin{array}{c}
A:H^{\kappa+s}_{\ell}(\mathbb{R}^{\kappa})\longrightarrow H^{\kappa+s}_{\ell}(\mathbb{R}^{\kappa})\\[0.15in]
{\displaystyle f\longmapsto
(Af)(x,\mu)=
\int_{\Gamma_0}\,\gamma\,f(\gamma^{-1}x,\mu)\,d\gamma}
\end{array}
\]
\label{Adef}
\end{Def}
\begin{prop}
$A$ is a projection.  Furthermore,
\begin{equation}
\mbox{\rm range}\,A=
H^{\kappa+s}_{\ell}(\mathbb{R}^{\kappa},\mathbb{T}^p)
\label{Arange}
\end{equation}
and
\begin{equation}
\mbox{\rm ker}\,A=
\mbox{\rm range}\,{\cal L}_{B}
\label{Aker}
\end{equation}
\label{prop_Adecomposition}
\end{prop}
\proof
The proof is given in the appendix.

Since $A$ is a projection, then
$H^{\kappa+s}_{\ell}(\mathbb{R}^{\kappa})=\mbox{\rm
  ker}\,A\oplus\mbox{\rm range}\,A$, and Proposition
\ref{prop_Adecomposition} shows that this decomposition is precisely
the decomposition of $H^{\kappa+s}_{\ell}(\mathbb{R}^{\kappa})$ given
in Proposition \ref{prop_enf2}.  Thus, for any $f\in
H^{\kappa+s}_{\ell}(\mathbb{R}^{\kappa})$, write
\[
f=Af+(I-A)f,
\]
and note that $Af$ is $\mathbb{T}^p$-equivariant and that $(I-A)f\in
\mbox{\rm ker}\,A$.  
From Proposition \ref{prop_Adecomposition}, there exists
$h\in H^{\kappa+s}_{\ell}(\mathbb{R}^{\kappa})$ such that ${\cal
  L}_{B}h=(I-A)f$.

\subsection{Phase decoupling}
The following example serves as an illustration of a trivial
(well-known) case in
which normal form toroidal symmetry leads to a decoupling of the equations in
the normal form.
\begin{examp}
In the case where $B=\mbox{\rm diag}(i\omega,-i\omega)$, the normal
form has the rotational symmetry of a one-dimensional torus:
$(x,\overline{x})\rightarrow\,(e^{i\theta}\,x,e^{-i\theta}\overline{x})$,
$\theta\in\mathbb{T}^1$, and it is
easy to verify that the most general ${\mathbb T}^1$-equivariant
differential equation has the form
\begin{equation}
\dot{x}=f(x\overline{x})x,\,\,\,\,(x\in\mathbb{C})
\label{t1nf}
\end{equation}
and its complex conjugate
where $f$ is complex-valued, and $f(0)=i\omega$.  Writing
$x=re^{i\theta}$ leads to the equations
\begin{equation}
\dot{r}=\mbox{\rm
  Re}(f(r^2))\,r,\,\,\,\,\,\,\,\,\,\dot{\theta}=\mbox{\rm Im}(f(r^2)).
\label{t1nfuc}
\end{equation}
We note that $\theta$ does not appear in the $\dot{r}$ equation, and
that the $\dot{r}$ equation has a reflectional symmetry $r\rightarrow
-r$.  The analysis of the normal form (\ref{t1nf}) then essentially reduces to a
one-dimensional problem (the $\dot{r}$ equation in (\ref{t1nfuc})) which possesses some
residual (reflectional) symmetry.
\end{examp}

In fact, this example is a special case of a more general result which
holds when the spectrum of $B$ satisfies Hypothesis \ref{spectralhyp}(b), and which we now
outline.

From Lemma \ref{lemenf1}, we get the following
\begin{cor}
Suppose the spectrum of $B$ satisfies Hypothesis \ref{spectralhyp}(b).  Then a smooth vector field
$f:\mathbb{R}^{\kappa+s}\longrightarrow\mathbb{R}^{\kappa}$
with $f(0,0)=0$, $Df(0,0)=0$
is $\Gamma_0\cong \mathbb{T}^p$-equivariant if and only if
$f$ has the form
\begin{equation}
f(x,\mu)=\left\{\begin{array}{lcl}
\left(
\begin{array}{c}
a_1(x_1\overline{x_1},\ldots,x_p\overline{x_p},\mu)\,x_1\\ \\
\overline{a_1(x_1\overline{x_1},\ldots,x_p\overline{x_p},\mu)\,x_1}\\ \\
\vdots\\ \\
a_p(x_1\overline{x_1},\ldots,x_p\overline{x_p},\mu)\,x_p\\ \\
\overline{a_p(x_1\overline{x_1},\ldots,x_p\overline{x_p},\mu)\,x_p},
\end{array}
\right)&\mbox{\rm if}&\kappa=2p\\
&&\\
\left(
\begin{array}{c}
a_0(x_0,x_1\overline{x_1},\ldots,x_p\overline{x_p},\mu)\\ \\
a_1(x_0,x_1\overline{x_1},\ldots,x_p\overline{x_p},\mu)\,x_1\\ \\
\overline{a_1(x_0,x_1\overline{x_1},\ldots,x_p\overline{x_p},\mu)\,x_1}\\ \\
\vdots\\ \\
a_p(x_0,x_1\overline{x_1},\ldots,x_p\overline{x_p},\mu)\,x_p\\ \\
\overline{a_p(x_0,x_1\overline{x_1},\ldots,x_p\overline{x_p},\mu)\,x_p},
\end{array}
\right)&\mbox{\rm if}&\kappa=2p+1,\end{array}\right.
\label{torus_nf2}
\end{equation}
where $a_1,\ldots,a_p$ are smooth and complex-valued, and $a_{0}$ is
smooth and real-valued.  
\end{cor}
\begin{prop}
Consider a differential equation
$\dot{x}=Bx+f(x,\mu)$, where $f$ is as in
(\ref{torus_nf2}).  Then under the under the change
of variables $x_0=\rho_0$, $x_j=\rho_je^{i\theta_j}$, $j=1,\ldots,p$, this
differential equation transforms into
\begin{equation}
\begin{array}{lcl}
\dot{\rho}_j=\mbox{\rm
  Re}(a_j(\rho_1^2,\ldots,\rho_p^2,\mu))\,\rho_j,\,\,\,j=1,\ldots,p&\mbox{\rm if}&\kappa=2p\\
&&\\
\left\{\begin{array}{l}
\dot{\rho}_0=a_0(\rho_0,\rho_1^2,\ldots,\rho_p^2,\mu)\\ \\
\dot{\rho}_j=\mbox{\rm
  Re}(a_j(\rho_0,\rho_1^2,\ldots,\rho_p^2,\mu))\,\rho_j,\,\,\,j=1,\ldots,p
\end{array}\right.&\mbox{\rm if}&\kappa=2p+1,\end{array}
\label{radial_eqs}
\end{equation}
and
\begin{equation}
\dot{\theta}_j=\left\{
\begin{array}{lcl}
\mbox{\rm
  Im}(a_j(\rho_1^2,\ldots,\rho_p^2,\mu)),\,\,\,j=1,\ldots,p&\mbox{\rm
  if}&\kappa=2p\\
&&\\
\mbox{\rm
  Im}(a_j(\rho_0,\rho_1^2,\ldots,\rho_p^2,\mu)),\,\,\,j=1,\ldots,p&\mbox{\rm if}&\kappa=2p+1.\end{array}\right.
\label{angular_eqs}
\end{equation}
\end{prop}
\proof
This is a simple computation.
\hfill\qed

We will call the subsystem (\ref{radial_eqs}) the {\em uncoupled
  radial part} of the normal form (\ref{torus_nf2}).
For many practical purposes of interest, it is
sufficient to consider only the uncoupled radial part
  (\ref{radial_eqs}) in the analysis of (\ref{torus_nf2}).
For example, small-amplitude equilibria of (\ref{radial_eqs})
  correspond to periodic solutions or invariant
tori of the full normal form (\ref{torus_nf2}).  Oftentimes, given some
  normal hyperbolicity conditions, these
invariant objects for (\ref{torus_nf2}) persist as invariant objects in the original system
(\ref{y7}).  In fact, in the case
of non-resonant double Hopf bifurcation
($\Lambda_0=\{\pm\,i\omega_1,\pm\,i\omega_2\}$) and in the case of
saddle-node/Hopf interaction ($\Lambda_0=\{0,\pm\,i\omega_1\}$), it is
well-known \cite{Ta} that given some generic non-degeneracy conditions on the
coefficients of the lower-order nonlinearities, the radial equations
(\ref{radial_eqs}) (suitably truncated) completely determine the
dynamics in the full system (\ref{y7}) up to topological equivalence.
So, it is reasonable to investigate the
realizability of the uncoupled radial part (\ref{radial_eqs}) by center manifold reduction
(\ref{y6}) of
the RFDE (\ref{y1}).  

We now introduce an integer $d$ which should be interpreted such that $d=p$ in
the case where $\Lambda_0=\{\pm i\omega_1,\ldots,\pm i\omega_p\}$, and
$d=p+1$ in the case where $\Lambda_0=\{0,\pm i\omega_1,\ldots,\pm i\omega_p\}$.
Denote by ${\mathbb Z}_{2,p}$ the group
whose action on $\mathbb{R}^{d}$ is given by
\begin{equation}
\begin{array}{lcl}
(\rho_1,\ldots,\rho_p)\rightarrow(\lambda_1\rho_1,\ldots,\lambda_p\rho_p)&\mbox{\rm
  if}&d=p\\
&&\\
(\rho_0,\rho_1,\ldots,\rho_p)\rightarrow(\rho_0,\lambda_1\rho_1,\ldots,\lambda_p\rho_p)&\mbox{\rm if}&d=p+1,\end{array}
\label{multiz2act}
\end{equation}
where $\lambda_j\in\{1,-1\}$, $j=1,\ldots,p$.
\begin{Def}
For a given integer $\ell\geq 2$, a given normed space $X$, and for $d=p$
(respectively $d=p+1$), we denote
by $H^{d+s}_{\ell}(X)$ the linear space of homogeneous polynomials
of degree $\ell$ in the $d+s$ variables $\rho=(\rho_1,\ldots,\rho_p)$ (respectively
$\rho=(\rho_0,\rho_1,\ldots,\rho_p)$ and
$\mu=(\mu_1,\ldots,\mu_s)$ with coefficients in $X$. 
Denote by
$H^{d+s}_{\ell}(\mathbb{R}^{d},\mathbb{Z}_{2,p})\subset H^{d+s}_{\ell}(\mathbb{R}^{d})$
the subspace of $H^{d+s}_{\ell}(\mathbb{R}^{d})$ 
consisting of
$\mathbb{Z}_{2,p}$-equivariant polynomials.
\end{Def}
It is easy to show (see \cite{GSSII}) that the most general element of
$H^{d+s}_{\ell}(\mathbb{R}^{d},\mathbb{Z}_{2,p})$ has the form
\[
\begin{array}{lcl}
\left(
\begin{array}{c}
h_1(\rho_1^2,\ldots,\rho_p^2,\mu)\,\rho_1\\
\vdots\\
h_p(\rho_1^2,\ldots,\rho_p^2,\mu)\,\rho_p,
\end{array}
\right)&\mbox{\rm if}&d=p\\
&&\\
\left(
\begin{array}{c}
h_0(\rho_0,\rho_1^2,\ldots,\rho_p^2,\mu)\\
h_1(\rho_0,\rho_1^2,\ldots,\rho_p^2,\mu)\,\rho_1\\
\vdots\\
h_p(\rho_0,\rho_1^2,\ldots,\rho_p^2,\mu)\,\rho_p,
\end{array}
\right)&\mbox{\rm if}&d=p+1\end{array}
\]
and one immediately notices the similarity with (\ref{radial_eqs}).
It then becomes useful to define the following surjective linear mapping
\begin{equation}
\Pi : H^{\kappa+s}_{\ell}(\mathbb{R}^{\kappa},\mathbb{T}^p)\longrightarrow
H^{d+s}_{\ell}(\mathbb{R}^{d},\mathbb{Z}_{2,p})
\label{Pidef1}
\end{equation}
which is defined by sending the general element (\ref{torus_nf2}) of $H^{\kappa+s}_{\ell}(\mathbb{R}^{\kappa},\mathbb{T}^p)$
to the following element of $H^{d+s}_{\ell}(\mathbb{R}^{d},\mathbb{Z}_{2,p})$:
\begin{equation}
\begin{array}{lcl}
\left(
\begin{array}{c}
\mbox{\rm Re}(a_1(\rho_1^2,\ldots,\rho_p^2,\mu))\,\rho_1\\
\vdots\\
\mbox{\rm Re}(a_p(\rho_1^2,\ldots,\rho_p^2,\mu))\,\rho_p
\end{array}
\right)&\mbox{\rm if}&d=p\\
&&\\
\left(
\begin{array}{c}
a_0(\rho_0,\rho_1^2,\ldots,\rho_p^2,\mu)\\
\mbox{\rm Re}(a_1(\rho_0,\rho_1^2,\ldots,\rho_p^2,\mu))\,\rho_1\\
\vdots\\
\mbox{\rm Re}(a_p(\rho_0,\rho_1^2,\ldots,\rho_p^2,\mu))\,\rho_p
\end{array}
\right)&\mbox{\rm if}&d=p+1.\end{array}
\label{Pidef2}
\end{equation}

The following characterization of the mapping $\Pi$ will be
very useful later for computational purposes:
\mbox{\rm if $G$} is
an element of $H^{\kappa+s}_{\ell}(\mathbb{R}^{\kappa},\mathbb{T}^p)$, then 
\begin{equation}
(\Pi G)(\rho,\mu)=
C\cdot\,\gamma\,\cdot\,G(\gamma^{-1}\cdot R,\mu),
\label{Pidef3}
\end{equation}
where $\gamma$ is any fixed element of $\Gamma_0$,
$R=(\rho_1,\rho_1,\ldots,\rho_p,\rho_p)$ if $d=p$ and
$R=(\rho_0,\rho_1,\rho_1,\ldots,\rho_p,\rho_p)$ if $d=p+1$, and
where $C$ is the following $d\times \kappa$ matrix
\begin{equation}
C=\left\{\begin{array}{lcl}
\left(\begin{array}{cccccccc}
1/2&1/2&0&0&\cdots&\cdots&0&0\\
0&0&1/2&1/2&\cdots&\cdots&0&0\\
&&\vdots&\vdots&&&&\\
0&0&0&0&\cdots&\cdots&1/2&1/2
\end{array}
\right)&\mbox{\rm if}&d=p,\kappa=2p\\
&&\\
\left(\begin{array}{ccccccccc}
1&0&0&0&0&\cdots&\cdots&0&0\\
0&1/2&1/2&0&0&\cdots&\cdots&0&0\\
0&0&0&1/2&1/2&\cdots&\cdots&0&0\\
&&&\vdots&\vdots&&&&\\
0&0&0&0&0&\cdots&\cdots&1/2&1/2
\end{array}
\right)&\mbox{\rm if}&d=p+1,\kappa=2p+1.\end{array}\right.
\label{Cdef}
\end{equation}

\subsection{Parameter splitting}

There is a canonical direct sum decomposition of $H^{\kappa+s}_{\ell}(X)$
which will turn out to be quite useful for our purposes.  Note that
$H^{\kappa+s}_{\ell}(X)$ contains $H^{\kappa}_{\ell}(X)$ (the
$\mu$-independent polynomials) as a subspace, and consequently we can
write
\begin{equation}
H^{\kappa+s}_{\ell}(X)=H^{\kappa}_{\ell}(X)\oplus P^{\kappa+s}_{\ell}(X),
\label{paramnoparam}
\end{equation}
where $q\in P^{\kappa+s}_{\ell}(X)$ if and only if $q\in
H^{\kappa+s}_{\ell}(X)$ and $q(x,0)=0$.

The homological operator ${\cal L}_B$ (see (\ref{homdefine2})) preserves the decomposition (\ref{paramnoparam}):
\[
{\cal L}_B(H^{\kappa}_{\ell}(\mathbb{R}^{\kappa}))\subset
  H^{\kappa}_{\ell}(\mathbb{R}^{\kappa})\,\,\,\,\,\mbox{\rm and}\,\,\,\,\,
{\cal L}_B(P^{\kappa+s}_{\ell}(\mathbb{R}^{\kappa}))\subset
P^{\kappa+s}_{\ell}(\mathbb{R}^{\kappa}).
\]
Moreover,
\begin{equation}
H^{\kappa+s}_{\ell}(\mathbb{R}^{\kappa},\mathbb{T}^p)=H^{\kappa}_{\ell}(\mathbb{R}^{\kappa},\mathbb{T}^p)\oplus
P^{\kappa+s}_{\ell}(\mathbb{R}^{\kappa},\mathbb{T}^p),
\label{y15}
\end{equation}
where
$H^{\kappa}_{\ell}(\mathbb{R}^{\kappa},\mathbb{T}^p)=H^{\kappa}_{\ell}(\mathbb{R}^{\kappa})\cap
H^{\kappa+s}_{\ell}(\mathbb{R}^{\kappa},\mathbb{T}^p)$ and
$P^{\kappa+s}_{\ell}(\mathbb{R}^{\kappa},\mathbb{T}^p)=P^{\kappa+s}_{\ell}(\mathbb{R}^{\kappa})\cap H^{\kappa+s}_{\ell}(\mathbb{R}^{\kappa},\mathbb{T}^p)$.

If ${\cal L}_B|_1$ and ${\cal L}_B|_2$ represent
respectively the restrictions of ${\cal L}_B$ on
$H^{\kappa}_{\ell}(\mathbb{R}^{\kappa})$ and on
$P^{\kappa+s}_{\ell}(\mathbb{R}^{\kappa})$, then we have
the following refinement of Proposition \ref{prop_enf2}:
\begin{prop}
\[
\begin{array}{l}
H^{\kappa}_{\ell}(\mathbb{R}^{\kappa})=H^{\kappa}_{\ell}(\mathbb{R}^{\kappa},\mathbb{T}^p)\oplus
\mbox{\rm range}\,{\cal L}_B|_1,\\ \\
P^{\kappa+s}_{\ell}(\mathbb{R}^{\kappa})=P^{\kappa+s}_{\ell}(\mathbb{R}^{\kappa},\mathbb{T}^p)\oplus
\mbox{\rm range}\,{\cal L}_B|_2.
\end{array}
\]
\label{prop_enf3}
\end{prop}

The equivariant projection operator $A$ defined in Definition
\ref{Adef} also preserves the decomposition (\ref{paramnoparam}),
and we get the following refinement of Proposition \ref{prop_Adecomposition}
\begin{prop}
\begin{equation}
A(H^{\kappa}_{\ell}(\mathbb{R}^{\kappa}))=H^{\kappa}_{\ell}(\mathbb{R}^{\kappa},\mathbb{T}^p),\,\,\,\,\,\,\,
A(P^{\kappa+s}_{\ell}(\mathbb{R}^{\kappa}))=P^{\kappa+s}_{\ell}(\mathbb{R}^{\kappa},\mathbb{T}^p).
\label{Adsdecom}
\end{equation}
\label{prop_Aproj}
If $A|_1$ and $A|_2$ represent
respectively the restrictions of $A$ on
$H^{\kappa}_{\ell}(\mathbb{R}^{\kappa})$ and on
$P^{\kappa+s}_{\ell}(\mathbb{R}^{\kappa})$, then
\begin{equation}
\begin{array}{c}
\mbox{\rm ker}\,A|_1=\mbox{\rm range}\,{\cal L}_B|_1\\ \\
\mbox{\rm ker}\,A|_2=\mbox{\rm range}\,{\cal L}_B|_2
\end{array}
\label{y12}
\end{equation}
\label{prop_dual}
\end{prop}
\begin{rmk}
We note that
there exist similar direct sum decompositions of
$H^{d+s}_{\ell}(\mathbb{R}^{d})$
and of $H^{d+s}_{\ell}(\mathbb{R}^{d},\mathbb{Z}_{2,p})$ using the
subspace $H^{d}_{\ell}(\mathbb{R}^{d})\subset H^{d+s}_{\ell}(\mathbb{R}^{d})$
of $\mu$-independent polynomials
\begin{equation}
\begin{array}{c}
H^{d+s}_{\ell}(\mathbb{R}^{d})=H^{d}_{\ell}(\mathbb{R}^{d})\oplus
P^{d+s}_{\ell}(\mathbb{R}^{d})\\[0.15in]
H^{d+s}_{\ell}(\mathbb{R}^{d},\mathbb{Z}_{2,p})=
H^{d}_{\ell}(\mathbb{R}^{d},\mathbb{Z}_{2,p})\oplus
P^{d+s}_{\ell}(\mathbb{R}^{d},\mathbb{Z}_{2,p})
\end{array}
\label{y16}
\end{equation}
where $q\in P^{d+s}_{\ell}(\mathbb{R}^{d})$ if and only if
$q\in H^{d+s}_{\ell}(\mathbb{R}^{d})$ and $q(x,0)=0$, and where
$H^{d}_{\ell}(\mathbb{R}^{d},\mathbb{Z}_{2,p})=H^{d}_{\ell}(\mathbb{R}^{d})\cap
H^{d+s}_{\ell}(\mathbb{R}^{d},\mathbb{Z}_{2,p})$ and
$P^{d+s}_{\ell}(\mathbb{R}^{d},\mathbb{Z}_{2,p})=P^{d+s}_{\ell}(\mathbb{R}^{d})\cap
H^{d+s}_{\ell}(\mathbb{R}^{d},\mathbb{Z}_{2,p})$.
\label{dsdecomH2}
Note that the mapping $\Pi$ defined in (\ref{Pidef1})-(\ref{Pidef3}) 
preserves (\ref{y15}) and (\ref{y16}), i.e.
\[
\Pi(H^{\kappa}_{\ell}(\mathbb{R}^{\kappa},\mathbb{T}^p))=H^{d}_{\ell}(\mathbb{R}^{d},\mathbb{Z}_{2,p})\,\,\,\,\,\mbox{\rm
  and}\,\,\,\,\,
\Pi(P^{\kappa+s}_{\ell}(\mathbb{R}^{\kappa},\mathbb{T}^p))=P^{d+s}_{\ell}(\mathbb{R}^{d},\mathbb{Z}_{2,p}).
\]
\end{rmk}

Combining the results of this section with the Faria and
Magalh$\tilde{\mbox{\rm a}}$es normal form procedure described in
section 2, we get the following version of Theorem 5.8 of \cite{FMTB}
and Theorem 2.16 of \cite{FMH} which is adapted for our purposes
\begin{thm}
Consider the system (\ref{y4})
\begin{equation}
\begin{array}{rcl}
\dot{x}&=&{\dps Bx+\sum_{j\geq 2}\,f_j^1(x,y,\mu)}\\&&\\
{\dps\frac{d}{dt}y}&=&{\dps A_{Q^1}+\sum_{j\geq 2}\,f_j^2(x,y,\mu).}
\end{array}
\label{y7tayl}
\end{equation}
Write 
\begin{equation}
f_j^1(x,0,\mu)=h_j(x)+q_j(x,\mu), 
\label{fy0decom}
\end{equation}
where
$h_j\in H^{\kappa}_j(\mathbb{R}^{\kappa})$ and $q_j\in P^{\kappa+s}_j(\mathbb{R}^{\kappa})$.  
Then there is a formal near-identity change of variables
\[
(x,y)\longrightarrow
(\hat{x},\hat{y})+(U^1(\hat{x}),U^2(\hat{x}))+(W^1(\hat{x},\mu),W^2(\hat{x},\mu))
\]
(where $W^1(\hat{x},0)=0$, $W^2(\hat{x},0)=0$)
which transforms
(\ref{y7tayl}) into system (\ref{y5}) (upon dropping the hats),
and the flow on the invariant local center manifold $y=0$ is given by
\begin{equation}
\dot{x}=Bx+\sum_{j\geq 2}\,\left(\,(A|_1(h_j+Y_j))(x)+(A|_2(q_j+Z_j))(x,\mu)\right)
\label{y13}
\end{equation}
where $Y_2=0$, $Z_2=0$, and for $j\geq 3$, $Y_j(x)$ and $Z_j(x,\mu)$ are the
extra contributions to the terms of order $j$ coming from the
transformation of the lower order $(<j)$ terms, and $Z_j(x,0)=0$.
\label{thmdsdecomnf}
\end{thm}

\Section{Realizability: linear analysis}

In this section, we present the first of our main results on the realizability of
the radial part (\ref{Pidef2}) of toroidal normal forms
(\ref{torus_nf2}) {\em to any order}
for RFDEs (\ref{y4}) via the center-manifold normal form equations
(\ref{y13}).  

We will define a linear operator between suitable spaces of
polynomials, which arises in the context of the normal form
transformations of (\ref{y4}).  Our main result in this section will
be to establish the surjectivity of this operator.  Surjectivity will
be the main ingredient in the proof of our main realizability results
which will be presented in the next section.

Again, we we will try as much as possible to use concise notation
which will allow for the simultaneous treatment of both cases for
$\Lambda_0$ in Hypothesis \ref{spectralhyp}(b).

\subsection{Preliminaries}

For given integers $p\geq 1$ and $\ell\geq 2$, and for $\kappa=2p$
(respectively $2p+1$) and $d=p$ (respectively $p+1$), 
recall that $H^{d+s}_{\ell}(\mathbb{R}^d,\mathbb{Z}_{2,p})$ is
the linear space of homogeneous polynomials
of degree $\ell$ in the $d+s$ variables
$\rho=(\rho_1,\ldots,\rho_p)$ (respectively
$\rho=(\rho_0,\rho_1,\ldots,\rho_p)\equiv (\rho_0,\tilde{\rho})$) and
$\mu=(\mu_1,\ldots,\mu_s)$ with coefficients in $\mathbb{R}^{d}$,
and which are equivariant with respect to
  the $\mathbb{Z}_{2,p}$ action (\ref{multiz2act}) on $\mathbb{R}^d$.
Recall also that $H^{\kappa+s}_{\ell}(\mathbb{R}^{\kappa})$
is the
space of homogeneous polynomials of degree $\ell$ 
in the $\kappa+s$ variables
$x=(x_1,\overline{x_1},\ldots,x_p,\overline{x_p})$ (respectively
$x=(x_0,x_1,\overline{x_1},\ldots,x_p,\overline{x_p})$) and $\mu=(\mu_1,\ldots,\mu_s)$
with coefficients in $\mathbb{R}^{\kappa}$, and
$H^{\kappa+s}_{\ell}(\mathbb{R}^{\kappa},\mathbb{T}^p)$ is the subset of $H^{\kappa+s}_{\ell}(\mathbb{R}^{\kappa})$
consisting of $\mathbb{T}^p$-equivariant mappings.
\begin{Def}
Denote by $V^{d+s}_{\ell}(\mathbb{R})\subset H^{d+s}_{\ell}(\mathbb{R})$ the subspace of homogeneous degree
$\ell$ polynomials in the $d+s$ variables $v=(v_1,\ldots,v_p)$
(respectively $v=(v_0,v_1,\ldots,v_p)\equiv (v_0,\tilde{v})$) and
$\mu=(\mu_1,\ldots,\mu_s)$ with real coefficients, spanned by the basis
\begin{equation}
\begin{array}{lcl}
\{\,\mu^q\,
v^{2k}\,v_{c}\,\,|\,\,c\in\,\{1,\ldots,p\},\,\,(k,q)\in\mathbb{N}_0^{d+s},|q|+2|k|+1=\ell\}&\mbox{\rm
  if}&d=p\\
&&\\
\{\,\mu^q\,v_0^{k_0}\,\tilde{v}^{2k}\,\tilde{v}_c\,\,|\,\,c\in\,\{1,\ldots,p\},\,\,((k_0,k),q)\in\mathbb{N}_0^{d+s},|q|+k_0+2|k|+1=\ell\}\,\bigcup&&\\
\{\,\mu^q\,v_0^{k_0}\,\tilde{v}^{2k}\,\,|\,\,((k_0,k),q)\in\mathbb{N}_0^{d+s},|q|+k_0+2|k|=\ell\}&\mbox{\rm
  if}&d=p+1,
\end{array}
\label{Bbasis}
\end{equation}
where it is understood that if $(k,q)=(k_1,\ldots,k_p,q_1,\ldots,q_s)\in\mathbb{N}_0^{p+s}$, then $\mu^q=(\mu_1)^{q_1}\cdots (\mu_s)^{q_s}$,
$v^{2k}=\tilde{v}^{2k}=(v_1)^{2k_1}\cdots (v_p)^{2k_p}$,
$|q|=\sum_j\,q_j$ and $|k|=\sum_j\,k_j$.
\label{Vspacedef}
\end{Def}
Note that, $V^{d+s}_{\ell}(\mathbb{R})$ is isomorphic to the vector space
$H^{d+s}_{\ell}(\mathbb{R}^d,\mathbb{Z}_{2,p})$, since this latter
space has the following basis
\begin{equation}
\begin{array}{lcl}
\{\,\mu^q\,\rho^{2k}\,\rho_{c}\,\mbox{\bf
  e}_{c}\,\,|\,\,c\in\,\{1,\ldots,p\},\,\,(k,q)\in\mathbb{N}_0^{d+s},|q|+2|k|+1=\ell\}&\mbox{\rm
  if}&d=p\\
&&\\
\{\,\mu^q\,\rho_0^{k_0}\,\tilde{\rho}^{2k}\,\tilde{\rho}_c\,\mbox{\bf
  e}_{c+1}\,\,|\,\,c\in\,\{1,\ldots,p\},\,\,((k_0,k),q)\in\mathbb{N}_0^{d+s},|q|+k_0+2|k|+1=\ell\}\,\bigcup&&\\
\{\,\mu^q\,\rho_0^{k_0}\,\tilde{\rho}^{2k}\,\mbox{\bf
  e}_1\,\,|\,\,((k_0,k),q)\in\mathbb{N}_0^{d+s},|q|+k_0+2|k|=\ell\}&\mbox{\rm
  if}&d=p+1,
\end{array}
\label{Brbasis}
\end{equation}
where
$\mbox{\bf e}_{j}$ is a column vector with zeros on each row
except the $j^{\mbox{\small th}}$ row, which is 1.  Therefore,
$\mbox{\rm
  dim}\,H^{d+s}_{\ell}(\mathbb{R}^d,\mathbb{Z}_{2,p})=\mbox{\rm dim}\,V^{d+s}_{\ell}(\mathbb{R})$.

Since $B$ is as in (\ref{BdiagHopf}), then this corresponds to the
following choice of basis for the center subspace $P$:
\[
\Phi(t)=\left\{
\begin{array}{lcl}
(e^{i\omega_1t}\,\,\,\,e^{-i\omega_1t}\,\,\,\,\cdots\,\,\,\,e^{i\omega_pt}\,\,\,\,e^{-i\omega_pt})&\mbox{\rm
  if}&d=p,\,\kappa=2p\\
&&\\
(1\,\,\,\,e^{i\omega_1t}\,\,\,\,e^{-i\omega_1t}\,\,\,\,\cdots\,\,\,\,e^{i\omega_pt}\,\,\,\,e^{-i\omega_pt})&\mbox{\rm
  if}&d=p+1,\,\kappa=2p+1.
\end{array}\right.
\]
It follows that $\Psi(0)$ in
(\ref{projfdep}) is a $\kappa\times 1$ matrix
\begin{equation}
\Psi(0)=\left\{\begin{array}{lcl}\mbox{\rm
    col}(u_1,\overline{u_1},\ldots,u_p,\overline{u_p})&\mbox{\rm
    if}&\kappa=2p\\
&&\\
\mbox{\rm
    col}(u_0,u_1,\overline{u_1},\ldots,u_p,\overline{u_p})&\mbox{\rm if}&\kappa=2p+1,\end{array}\right.
\label{Vdef}
\end{equation}
where $u_0\neq 0$ is real and $u_j\neq 0$ are complex, $j=1,\ldots,p$.

\subsection{Linear analysis}

\begin{Def}
Let ${\cal S}$ denote the normed real linear space of $d\times \kappa$ matrices of the
form
\begin{equation}
M=\left\{\begin{array}{lcl}\left(\begin{array}{ccccc}
\alpha_{1,1}&\overline{\alpha_{1,1}}&\cdots&\alpha_{1,p}&\overline{\alpha_{1,p}}\\
\alpha_{2,1}&\overline{\alpha_{2,1}}&\cdots&\alpha_{2,p}&\overline{\alpha_{2,p}}\\
\vdots&\vdots&\cdots&\vdots&\vdots\\
\alpha_{p,1}&\overline{\alpha_{p,1}}&\cdots&\alpha_{p,p}&\overline{\alpha_{p,p}}
\end{array}\right)&\mbox{\rm if}&d=p,\,\kappa=2p\\
&&\\
\left(\begin{array}{cccccc}
\alpha_{0,0}&\alpha_{0,1}&\overline{\alpha_{0,1}}&\cdots&\alpha_{0,p}&\overline{\alpha_{0,p}}\\
\alpha_{1,0}&\alpha_{1,1}&\overline{\alpha_{1,1}}&\cdots&\alpha_{1,p}&\overline{\alpha_{1,p}}\\
\vdots&\vdots&\vdots&\cdots&\vdots&\vdots\\
\alpha_{p,0}&\alpha_{p,1}&\overline{\alpha_{p,1}}&\cdots&\alpha_{p,p}&\overline{\alpha_{p,p}}
\end{array}\right)&\mbox{\rm if}&d=p+1,\,\kappa=2p+1,\end{array}\right.
\label{calSdef}
\end{equation}
(where the $\alpha_{i,0}$ are real and the $\alpha_{i,j}$, $j\geq 1$ are complex numbers) equipped with norm
$||M||=\mbox{\rm max}(|\alpha_{i,j}|)$.  
For any given $M\in
{\cal S}$, we define the
{\em $\ell$-mapping associated to $M$}
\[
{\cal J}_M^{\ell} : H^{d+s}_{\ell}(\mathbb{R})\longrightarrow H^{\kappa+s}_{\ell}(\mathbb{R}^{\kappa})
\]
by
\[
({\cal
  J}_M^{\ell}(h))(x,\mu)=\Psi(0)\,h\left(x\left(M\right)^T,\mu\right),
\]
where $x=(x_1,\overline{x_1},\ldots,x_p,\overline{x_p})$ (respectively
$x=(x_0,x_1,\overline{x_1},\ldots,x_p,\overline{x_p})$).
\label{amdef}
\end{Def}
  
Let $\tau=(\tau_1,\ldots,\tau_d)$
be a vector (as of yet unspecified) in $\mathbb{R}^d$.
Define
\begin{equation}
E_{\tau}=\left(\begin{array}{c}\Phi(\tau_1)\\\vdots\\\Phi(\tau_d)\end{array}\right)=\left\{\begin{array}{lcl}
\left(\begin{array}{ccccc}
e^{i\omega_1\tau_1}&e^{-i\omega_1\tau_1}&\cdots&e^{i\omega_p\tau_1}&e^{-i\omega_p\tau_1}\\
e^{i\omega_1\tau_2}&e^{-i\omega_1\tau_2}&\cdots&e^{i\omega_p\tau_2}&e^{-i\omega_p\tau_2}\\
\vdots&\vdots&\cdots&\vdots&\vdots\\
e^{i\omega_1\tau_p}&e^{-i\omega_1\tau_p}&\cdots&e^{i\omega_p\tau_p}&e^{-i\omega_p\tau_p}
\end{array}\right)&\mbox{\rm if}&d=p\\
&&\\
\left(\begin{array}{cccccc}
1&e^{i\omega_1\tau_1}&e^{-i\omega_1\tau_1}&\cdots&e^{i\omega_p\tau_1}&e^{-i\omega_p\tau_1}\\
1&e^{i\omega_1\tau_2}&e^{-i\omega_1\tau_2}&\cdots&e^{i\omega_p\tau_2}&e^{-i\omega_p\tau_2}\\
\vdots&\vdots&\vdots&\vdots&\vdots&\vdots\\
1&e^{i\omega_1\tau_{p+1}}&e^{-i\omega_1\tau_{p+1}}&\cdots&e^{i\omega_p\tau_{p+1}}&e^{-i\omega_p\tau_{p+1}}
\end{array}\right)&\mbox{\rm if}&d=p+1,\end{array}\right.
\label{Edefmulthopf}
\end{equation}
and note that $E_{\tau}$ belongs to the space ${\cal S}$ of Definition \ref{amdef}.
Define the linear mapping
\[
{\cal E}^{\ell}_{{\tau}}:H^{d+s}_{\ell}(\mathbb{R})\longrightarrow H^{\kappa+s}_{\ell}(\mathbb{R}^{\kappa})
\]
by
\begin{equation}
({\cal
  E}^{\ell}_{\tau}(h))(x,\mu)\equiv ({\cal J}_{E_{\tau}}^{\ell})(h)(x,\mu)
\label{Ecaldefmulthopf}
\end{equation}
where ${\cal J}_{E_{\tau}}^{\ell}$ is the $\ell$-mapping associated to $E_{\tau}$.

Now, let
$\Pi:H^{\kappa+s}_{\ell}(\mathbb{R}^{\kappa},\mathbb{T}^p)\longrightarrow H^{d+s}_{\ell}(\mathbb{R}^d,\mathbb{Z}_{2,p})$
be the mapping defined in (\ref{Pidef1})-(\ref{Pidef3}), and let
$A:H^{\kappa+s}_{\ell}(\mathbb{R}^{\kappa})\longrightarrow
H^{\kappa+s}_{\ell}(\mathbb{R}^{\kappa},\mathbb{T}^p)$ be the group averaging
operator defined in
Definition \ref{Adef}.  Our main result in this section is the following:
\begin{prop}
For an open and dense set ${\cal U}\subset\mathbb{R}^d$,
the following linear mapping is surjective for all $\tau\in {\cal U}$:
\[
\Pi\circ A\circ {\cal E}^{\ell}_{{\tau}} : H^{d+s}_{\ell}(\mathbb{R})\longrightarrow
H^{d+s}_{\ell}(\mathbb{R}^d,\mathbb{Z}_{2,p}).
\]
\label{MultipleHopfinv}
\end{prop}
\proof
Let $K$ be the $d\times d$ matrix such that $K_{j,k}$
is equal to $-1$ if $j+k>d+1$ and is equal to $1$
otherwise.  It is easy to row reduce $K$ to the identity matrix, so $K$ is
invertible.
Therefore, $K$ induces an automorphism of the space $H^{d+s}_{\ell}(\mathbb{R})$:
\begin{equation}
\begin{array}{c}
{\cal K}:H^{d+s}_{\ell}(\mathbb{R})\longrightarrow
H^{d+s}_{\ell}(\mathbb{R})\\\\
({\cal K}h)(v,\mu)=h(v\,(K)^T,\mu).
\end{array}
\label{automorphdef}
\end{equation}
Define
$\widehat{V}^{d+s}_{\ell}(\mathbb{R})\equiv
{\cal K}^{-1}(V^{d+s}_{\ell}(\mathbb{R}))$, where
$V^{d+s}_{\ell}(\mathbb{R})$ is as in Definition \ref{Vspacedef}, and let 
\[
{\cal
    N}^{\ell}_{\tau}:\widehat{V}^{d+s}_{\ell}(\mathbb{R})\longrightarrow H^{d+s}_{\ell}(\mathbb{R}^d,\mathbb{Z}_{2,p})
\]
be the restriction of 
$\Pi\circ A\circ {\cal E}^{\ell}_{{\tau}}$ to
$\widehat{V}^{d+s}_{\ell}(\mathbb{R})$.  Our approach to proving
Proposition \ref{MultipleHopfinv} will be to prove that there exists an open and dense set
of points ${\cal U}\subset\mathbb{R}^d$ such that ${\cal
  N}^{\ell}_{{\tau}}$ is invertible for all $\tau\in {\cal U}$.

If $\langle\,{\cal N}^{\ell}_{{\tau}}\,\rangle$ is any matrix representation of
${\cal N}^{\ell}_{{\tau}}$, then $\mbox{\rm
  det}(\langle\,{\cal N}^{\ell}_{{\tau}}\,\rangle)$ is a real-analytic
function of $\tau_1,\ldots,\tau_d$ (in fact, it is a polynomial in 
$\cos\,\omega_k\tau_q$ and $\sin\,\omega_k\tau_q$,
$k\in\{1,\ldots,p\},\,q\in\{1,\ldots,d\}$).  Therefore, if we can show that $\mbox{\rm
det}(\langle\,{\cal N}^{\ell}_{{\tau}}\,\rangle)$
is not identically zero, the conclusion is a trivial consequence of
this analyticity.  This amounts to showing that
there exists at least one point ${\tau}^*\in\mathbb{R}^d$ such that with $E_{{\tau}^*}$ as in (\ref{Edefmulthopf}), the
mapping ${\cal N}^{\ell}_{{\tau}^*}$ is invertible.  We will prove this last claim
with a sequence of five lemmas.  
\begin{lemma}
Let ${\cal S}$ be as in Definition \ref{amdef}.  If $M_*\in {\cal S}$ is such that the restriction of the map $\Pi\circ
A\circ {\cal J}_{M_*}^{\ell}$ to
$\widehat{V}^{d+s}_{\ell}(\mathbb{R})$:
\[
\Pi\circ A\circ
{\cal J}_{M_*}^{\ell}:\widehat{V}^{d+s}_{\ell}(\mathbb{R})\longrightarrow
H^{d+s}_{\ell}(\mathbb{R}^d,\mathbb{Z}_{2,p})
\]
is invertible, 
then there is a $\delta=\delta(M_*)>0$ such
that for all $M$ in the $\delta$-ball centered on $M_*$, the restriction
$\Pi\circ A\circ {\cal J}_M^{\ell}:\widehat{V}^{d+s}_{\ell}(\mathbb{R})\longrightarrow
H^{d+s}_{\ell}(\mathbb{R}^d,\mathbb{Z}_{2,p})$ is invertible.
\label{lemSdef}
\end{lemma}
\proof
This follows from the fact that the determinant of the map
$\Pi\circ A\circ {\cal J}^{\ell}_M:\widehat{V}^{d+s}_{\ell}(\mathbb{R})\longrightarrow
H^{d+s}_{\ell}(\mathbb{R}^d,\mathbb{Z}_{2,p})$ is continuous in the
entries of $M$.
\hfill\qed

As mentioned above, $\Psi(0)$ in (\ref{Vdef}) is such that each of its
components is non-zero.  Therefore:
\begin{lemma}
There exists $\sigma_1,\ldots,\sigma_p$ such that
\[
\mbox{\rm Re}(e^{i\sigma_j}u_j)\neq 0,\,\,\,j=1,\ldots,p.
\]
\label{lemphis}
\end{lemma} 
\begin{lemma}
Let ${\cal S}$ be as defined in Definition \ref{amdef} and
$\sigma_1,\ldots,\sigma_p$ be as in Lemma \ref{lemphis}.  Consider the
following element $I\in {\cal S}$ of the form (\ref{calSdef}) where
$\alpha_{0,0}=1$ in the case $d=p+1$, and:
\[
\alpha_{j,k}=\left\{\begin{array}{ccc}
e^{i\sigma_j}&\mbox{\rm if}&j=k\geq 1\\[0.15in]
0&\mbox{\rm if}&j\neq k
\end{array}\right.
\]
If ${\cal J}^{\ell}_I : H^{d+s}_{\ell}(\mathbb{R})\longrightarrow H^{\kappa+s}_{\ell}(\mathbb{R}^{\kappa})$
is the $\ell$-mapping associated to $I$, then the restriction to
$V^{d+s}_{\ell}(\mathbb{R})$:
\[
\mbox{\rm $\Pi\circ A\circ {\cal J}^{\ell}_I : V^{d+s}_{\ell}(\mathbb{R})\longrightarrow
H^{d+s}_{\ell}(\mathbb{R}^d,\mathbb{Z}_{2,p})$} 
\]
is invertible.
\label{lemI}
\end{lemma}
\proof
We give only the proof in the case $d=p$, $\kappa=2p$.  The other case
($d=p+1$, $\kappa=2p+1$) is treated in a completely similar manner.

Consider the basis element 
$\mu^q\,v^{2k}\,v_c$ of
$V^{p+s}_{\ell}(\mathbb{R})$
(see (\ref{Bbasis})).  Then after an appropriate translation of the
integration variables, we have
\[
\begin{array}{l}
(\Pi\circ A\circ {\cal J}^{\ell}_I)(\mu^q\,v^{2k}\,v_c)=\\[0.15in]
{\displaystyle
\frac{\mu^q}{(2\pi)^p}\int_0^{2\pi}\cdots\int_0^{2\pi}\,{\cal G}\, (\rho_1(e^{-i\theta_1}+e^{i\theta_1}))^{2k_1}\cdots
(\rho_p(e^{-i\theta_p}+e^{i\theta_p}))^{2k_p}(\rho_{c}(e^{-i\theta_{c}}+e^{i\theta_{c}}))d\theta_1\cdots
d\theta_p},
\end{array}
\]
where ${\cal G}=C\cdot\mbox{\rm
  diag}(e^{i\theta_1},e^{-i\theta_1},\ldots,e^{i\theta_p},e^{-i\theta_p})\cdot\mbox{\rm
  diag}(e^{i\sigma_1},e^{-i\sigma_1},\ldots,e^{i\sigma_p},e^{-i\sigma_p})\cdot 
  \Psi(0)$, ($C$ as
in (\ref{Cdef}) and $\Psi(0)$
  as in (\ref{Vdef})).
A simple computation then shows that
\[
\begin{array}{l}
(\Pi\circ A\circ {\cal J}^{\ell}_I)(\mu^q\,v^{2k}\,v_c)=\\[0.15in]
{\displaystyle \mbox{\rm
  Re}(e^{i\sigma_{c}}u_{c})\frac{2k_{c}+1}{k_{c}+1}\left[
\frac{(2k_1)!}{(k_1!)^2}\frac{(2k_2)!}{(k_2!)^2}\cdots\frac{(2k_p)!}{(k_p!)^2}\right]\,\mu^q\,\rho^{2k}\,\rho_{c}\,\mbox{\bf e}_{c}},
\end{array}
\]
where we remind the reader that
$\mbox{\bf e}_{c}$ is a $p$-dimensional column vector with zeros on each row
except the $c^{\mbox{\small th}}$ row, which is 1.
Taking into account (\ref{Bbasis}) and (\ref{Brbasis}), under a suitable choice of bases
for the spaces 
$V^{p+s}_{\ell}(\mathbb{R})$ and $H^{p+s}_{\ell}(\mathbb{R}^p,\mathbb{Z}_{2,p})$,
the matrix
  representation of the restriction of $\Pi\circ A\circ {\cal
    J}^{\ell}_I$ to
$V^{p+s}_{\ell}(\mathbb{R})$ is diagonal with non-zero diagonal
  entries.
\hfill\qed
\begin{lemma}
Let ${\cal S}$ be as in Definition \ref{amdef}, 
$\sigma_1,\ldots,\sigma_p$ as in Lemma \ref{lemphis} and $I$ as in
Lemma \ref{lemI}.  
Consider the
element $E_*\equiv KI\in {\cal S}$, where the $d\times d$ matrix $K$ is such that
$K_{j,k}$ is equal to $-1$ if $j+k>d+1$ and is equal to $1$
otherwise.
If \mbox{\rm ${\cal J}^{\ell}_{E_*} : H^{d+s}_{\ell}(\mathbb{R})\longrightarrow H^{\kappa+s}_{\ell}(\mathbb{R}^{\kappa})$}
is the $\ell$-mapping associated to $E_*$, then the restriction
\mbox{\rm $\Pi\circ A\circ {\cal J}^{\ell}_{E_*} : \widehat{V}^{d+s}_{\ell}(\mathbb{R})\longrightarrow
H^{d+s}_{\ell}(\mathbb{R}^d,\mathbb{Z}_{2,p})$} is
invertible.
\label{lemEstar}
\end{lemma}
\proof
Let $I$ and ${\cal J}^{\ell}_I$ be as in Lemma \ref{lemI}.  Since $E_*=KI$,
it follows that ${\cal J}^{\ell}_{E_*}={\cal J}^{\ell}_I\circ {\cal
  K}$, 
where ${\cal K}$ is the automorphism defined in (\ref{automorphdef}).
So 
$\Pi\circ A\circ {\cal J}^{\ell}_{E_*}=(\Pi\circ A\circ {\cal
  J}^{\ell}_I)\circ {\cal K}$.
Consequently, the restriction of $\Pi\circ A\circ {\cal
  J}^{\ell}_{E_*}$ to $\widehat{V}^{d+s}_{\ell}(\mathbb{R})\equiv {\cal K}^{-1}(V^{d+s}_{\ell}(\mathbb{R}))$ is invertible.
\hfill\qed
\begin{lemma}
Let $E_*$ be as in Lemma \ref{lemEstar} and let $\delta=\delta(E_*)>0$
be as in Lemma \ref{lemSdef}.  There exists a $\tau^*\in\mathbb{R}^p$ such that
$E_{{\tau}^*}$ in (\ref{Edefmulthopf}) satisfies
$||E_{{\tau}^*}-E_*||<\delta$, and consequently 
if ${\cal E}^{\ell}_{{\tau}^*}\equiv {\cal
  J}^{\ell}_{E_{{\tau}^*}} : H^{d+s}_{\ell}(\mathbb{R})\longrightarrow H^{\kappa+s}_{\ell}(\mathbb{R}^{\kappa})$
is the $\ell$-mapping associated to $E_{{\tau}^*}$,
then the restriction
${\cal N}^{\ell}_{{\tau}^*}\equiv\Pi\circ A\circ {\cal
  E}^{\ell}_{{\tau}^*}:\widehat{V}^{d+s}_{\ell}(\mathbb{R})\longrightarrow H^{d+s}_{\ell}(\mathbb{R}^d,\mathbb{Z}_{2,p})$ is invertible (from Lemma \ref{lemSdef}).
\label{finallem}
\end{lemma}
\proof
Since the $\omega_1,\ldots,\omega_p$ are independent over the
rationals, it follows that the
set
\[
\{(e^{i\omega_1 t},e^{i\omega_2 t},\ldots,e^{i\omega_p
  t})\,\,|\,\,t\in\mathbb{R}\}
\]
is dense on the $p$-torus $\mathbb{T}^p$.  Consequently, it is
possible to choose a $\tau^*\in\mathbb{R}^d$ such that each row of $E_{\tau^*}$ is as close (in any given norm) as we wish to the
corresponding row of $E_*$.
\hfill\qed

\noindent
The proof of Proposition \ref{MultipleHopfinv} follows immediately from
Lemmas \ref{lemSdef}-\ref{finallem}.
\hfill\qed

\subsection{Refinement}

In the next section, we will need a finer version of Proposition
\ref{MultipleHopfinv}.  If $H^{d}_{\ell}(\mathbb{R})$ denotes the
subspace of $\mu$-independent elements of 
$H^{d+s}_{\ell}(\mathbb{R})$, then we have
\begin{equation}
H^{d+s}_{\ell}(\mathbb{R})=H^{d}_{\ell}(\mathbb{R})\oplus P^{d+s}_{\ell}(\mathbb{R}),
\label{y20}
\end{equation}
where $q\in P^{d+s}_{\ell}(\mathbb{R})$ if and only if
$q\in H^{d+s}_{\ell}(\mathbb{R})$ and $q(v,0)=0$.
Define $V^d_{\ell}(\mathbb{R})=V^{d+s}_{\ell}(\mathbb{R})\cap
H^d_{\ell}(\mathbb{R})$ and
$W^d_{\ell}(\mathbb{R})=V^{d+s}_{\ell}(\mathbb{R})\cap
P^{d+s}_{\ell}(\mathbb{R})$ and note that
\begin{equation}
V^{d+s}_{\ell}(\mathbb{R})=V^d_{\ell}(\mathbb{R})\oplus
W^{d+s}_{\ell}(\mathbb{R})
\label{y22}
\end{equation}
is precisely the decomposition of $V^{d+s}_{\ell}(\mathbb{R})$ into
the direct sum of $\mu$-independent elements of
$V^{d+s}_{\ell}(\mathbb{R})$ and elements of
$V^{d+s}_{\ell}(\mathbb{R})$ which vanish at $\mu=0$.
The automorphism ${\cal K}$ of $H^{d+s}_{\ell}(\mathbb{R})$ defined in (\ref{automorphdef})
preserves these decompositions, and we have
\begin{equation}
\begin{array}{lll}
\widehat{V}^{d+s}_{\ell}(\mathbb{R})&=&{\cal
  K}^{-1}(V^{d+s}_{\ell}(\mathbb{R}))\\[0.15in]
&=&{\cal
  K}^{-1}(V^d_{\ell}(\mathbb{R}))\oplus {\cal
  K}^{-1}(W^{d+s}_{\ell}(\mathbb{R}))\\[0.15in]
&\equiv&\widehat{V}^d_{\ell}(\mathbb{R})\oplus\widehat{W}^{d+s}_{\ell}(\mathbb{R})
\end{array}
\label{y23}
\end{equation}
which is the decomposition of $\widehat{V}^{d+s}_{\ell}(\mathbb{R})$
into
the direct sum of $\mu$-independent elements of
$\widehat{V}^{d+s}_{\ell}(\mathbb{R})$ and elements of
$\widehat{V}^{d+s}_{\ell}(\mathbb{R})$ which vanish at $\mu=0$.

Then,
taking into account (\ref{y16}), we have
\begin{prop}
\begin{equation}
\mbox{\rm dim}\,\widehat{V}^d_{\ell}(\mathbb{R})=\mbox{\rm
  dim}\,H^{d}_{\ell}(\mathbb{R}^d,\mathbb{Z}_{2,p})\,\,\,\,\,\mbox{\rm and}\,\,\,\,\,
\mbox{\rm dim}\,\widehat{W}^{d+s}_{\ell}(\mathbb{R})=\mbox{\rm
dim}\,P^{d+s}_{\ell}(\mathbb{R}^{d},\mathbb{Z}_{2,p}).
\label{y21}
\end{equation}
Furthermore,
if $\tau$ is as in Proposition \ref{MultipleHopfinv}, then
\[
(\Pi\circ A\circ {\cal
  E}^{\ell}_{\tau})(\widehat{V}^d_{\ell}(\mathbb{R}))=H^d_{\ell}(\mathbb{R}^d,\mathbb{Z}_{2,p})\,\,\,\,\,\mbox{\rm
  and}\,\,\,\,\,
(\Pi\circ A\circ {\cal
  E}^{\ell}_{\tau})(\widehat{W}^{d+s}_{\ell}(\mathbb{R}))=P^{d+s}_{\ell}(\mathbb{R}^d,\mathbb{Z}_{2,p}).
\]
\label{dsdecomMultipleHopfinv}
\end{prop}
\proof
We give the proof in the case $d=p$, the other case being treated in
a similar manner.

We note that
\[
V^p_{\ell}(\mathbb{R})=\mbox{\rm span}\,\{\,
v^{2k}\,v_{c}\,\,|\,\,c\in\,\{1,\ldots,p\},\,\,k\in\mathbb{N}_0^{p},2|k|+1=\ell\},
\]
and
\[
H^{p}_{\ell}(\mathbb{R}^p,\mathbb{Z}_{2,p})=\mbox{\rm
  span}\,\{\,\rho^{2k}\,\rho_{c}\,\mbox{\bf
  e}_{c}\,\,|\,\,c\in\,\{1,\ldots,p\},\,\,k\in\mathbb{N}_0^{p},2|k|+1=\ell\}.
\]
Equation (\ref{y21}) follows from (\ref{y22}), (\ref{y23}) and the fact that $\widehat{V}^p_{\ell}(\mathbb{R})={\cal
  K}^{-1}(V^p_{\ell}(\mathbb{R}))$.
It now follows from the theory presented in section 3 that 
\[
(\Pi\circ A\circ {\cal
  E}^{\ell}_{\tau})(\widehat{V}^p_{\ell}(\mathbb{R}))\subset H^p_{\ell}(\mathbb{R}^p,\mathbb{Z}_{2,p})\,\,\,\,\,\mbox{\rm
  and}\,\,\,\,\,
(\Pi\circ A\circ {\cal
  E}^{\ell}_{\tau})(\widehat{W}^{p+s}_{\ell}(\mathbb{R}))\subset P^{p+s}_{\ell}(\mathbb{R}^p,\mathbb{Z}_{2,p}).
\]
The reverse inclusions then follow from the invertibility of ${\cal
  N}^{\ell}_{\tau}$ and from (\ref{y16}), (\ref{y20}) and (\ref{y21}).
\hfill\qed

\Section{Main Results}

We are now ready to state and prove our main realizability results for
both cases of Hypothesis \ref{spectralhyp}(b), with the convention
that respectively
$(d,\kappa)=(p,2p)$ and $(d,\kappa)=(p+1,2p+1)$.  It will be convenient to
define the following linear spaces of (non-homogeneous) polynomials
\begin{Def}
For an integer $\ell\geq 2$, define
\[
\begin{array}{ll}
\widehat{{\cal
  V}}^{d+s}_{\ell}(\mathbb{R})\equiv\oplus_{j=2}^{\ell}\,\widehat{V}^{d+s}_{j}(\mathbb{R}),\hspace*{1in}&
\widehat{{\cal
  V}}^d_{\ell}(\mathbb{R})\equiv\oplus_{j=2}^{\ell}\,\widehat{V}^d_j(\mathbb{R}),\\[0.15in]
\widehat{{\cal W}}^{d+s}_{\ell}(\mathbb{R})\equiv\oplus_{j=2}^{\ell}\,\widehat{W}^{d+s}_j(\mathbb{R}),&
 {\cal
  H}^{d+s}_{\ell}(\mathbb{R}^d,\mathbb{Z}_{2,p})\equiv\oplus_{j=2}^{\ell}\,H^{d+s}_{j}(\mathbb{R}^d,\mathbb{Z}_{2,p}),\\[0.15in]
{\cal
  H}^d_{\ell}(\mathbb{R}^d,\mathbb{Z}_{2,p})\equiv\oplus_{j=2}^{\ell}\,H^d_j(\mathbb{R}^d,\mathbb{Z}_{2,p}),&
{\cal
  P}^{d+s}_{\ell}(\mathbb{R}^d,\mathbb{Z}_{2,p})\equiv\oplus_{j=2}^{\ell}\,P^{d+s}_{j}(\mathbb{R}^d,\mathbb{Z}_{2,p}),\\[0.15in]
{\cal
  H}^d_{\ell}(\mathbb{R})\equiv\oplus_{j=2}^{\ell}\,H^d_j(\mathbb{R}),&
{\cal
  H}^{d+s}_{\ell}(\mathbb{R})\equiv\oplus_{j=2}^{\ell}\,H^{d+s}_{j}(\mathbb{R}).
\end{array}
\]
\label{nonhompolys}
\end{Def}

Our first result addresses the issue of realizability of singularities
and unfoldings within the class
of scalar delay-differential equations with $d$ delays.
\begin{thm}
Consider the RFDE (\ref{y1}), and let $\Lambda_0$ denote the set of
solutions of (\ref{ychar}) with zero real part.  Suppose that
Hypothesis \ref{spectralhyp} is satisfied. 
Let $\ell\geq 2$ be a given integer.  For each
$h\in {\cal H}^{d}_{\ell}(\mathbb{R}^d,\mathbb{Z}_{2,p})$:
\[
h(\rho)=\sum_{j=2}^{\ell}\,h_{j}(\rho),
\]
($h_{j}\in
H^{d}_{j}(\mathbb{R}^d,\mathbb{Z}_{2,p}),\,j=2,\ldots,\ell)$ and each
$q\in {\cal P}^{d+s}_{\ell}(\mathbb{R}^d,\mathbb{Z}_{2,p})$:
\[
q(\rho,\mu)=\sum_{j=2}^{\ell}\,q_j(\rho,\mu),
\]
($q_j\in P^{d+s}_j(\mathbb{R}^d,\mathbb{Z}_{2,p}),\,j=2,\ldots,\ell)$,
there are $d$ distinct points
$\tau_1,\ldots,\tau_d\in [-r,0]$, an $\eta\in \widehat{{\cal V}}^{d}_{\ell}(\mathbb{R})$:
\begin{equation}
\eta(v)=\sum_{j=2}^{\ell}\,\eta_{j}(v),
\label{etaform}
\end{equation}
($\eta_{j}\in \widehat{V}^{d}_j(\mathbb{R}),\,j=2,\ldots,\ell)$,
and a $\xi\in \widehat{{\cal W}}^{d+s}_{\ell}(\mathbb{R})$:
\begin{equation}
\xi(v,\mu)=\sum_{j=2}^{\ell}\,\xi_j(v,\mu),
\label{xiform}
\end{equation}
($\xi_j\in \widehat{W}^{d+s}_j(\mathbb{R}),\,j=2,\ldots,\ell)$,
such that if 
\[
\widetilde{F}(z_t,\mu)=\eta(z(t+\tau_1),\ldots,z(t+\tau_d))+\xi(z(t+\tau_1),\ldots,z(t+\tau_d),\mu)
\]
in
  (\ref{y1p}), then in polar coordinates, the radial part of the center manifold equations
  (\ref{y13}) in $\mathbb{T}^p$-equivariant normal form up to degree $\ell$ reduces to
  $\dot{\rho}=h(\rho)+q(\rho,\mu)$, where
  $\rho\equiv (\rho_1,\ldots,\rho_p)$ or $\rho=(\rho_0,\rho_1,\ldots,\rho_p)$.
  In fact, $\tau$ can be chosen in an open and
  dense set of $[-r,0]^d$, independently of the particular $h$ and $q$
  to be
  realized (i.e. only $\eta$ and $\xi$ must be changed in order to account for
  different jets to be realized).
\label{mainthm1}
\end{thm}
\proof
Choose a point ${\tau}\in [-r,0]^d$ such that the previously defined
linear mappings
\[
{\cal N}^{j}_{{\tau}}:\widehat{V}^{d+s}_{j}(\mathbb{R})\longrightarrow
H^{d+s}_{j}(\mathbb{R}^d,\mathbb{Z}_{2,p})
\]
are invertible for all $j=2,\ldots,\ell$ (from Proposition
\ref{MultipleHopfinv}, this is possible for an open and dense set of
points in $[-r,0]^d$).  Suppose $\eta$ is an arbitrary polynomial of
the form (\ref{etaform}), $\xi$ is and arbitrary polynomial of the
form (\ref{xiform}), and suppose $\widetilde{F}$ in (\ref{y1p}) is such that
$\widetilde{F}(z_t,\mu)=\eta(z(t+\tau_1),\ldots,z(t+\tau_d))+\xi(z(t+\tau_1),\ldots,z(t+\tau_d),\mu)$.
Using 
Theorem \ref{thmdsdecomnf}, it is
possible to define successively at each order near identity changes of
variables of the form 
\begin{equation}
(x,y)=(\hat{x},\hat{y})+(U_j^1(\hat{x})+W_j^1(\hat{x},\mu),U_j^2(\hat{x})+W_j^2(\hat{x},\mu)),
\label{dsdecomnfcv}
\end{equation}
where $W^i_j(\hat{x},0)=0,\,i=1,2$,
which transform (\ref{y4}) into
(\ref{y5}), and the center manifold equations are as in (\ref{y6}), 
with 
\begin{equation}
\begin{array}{rcl}
g_2^1(x,0,\mu)&=&A\, ({\cal E}^2_{\tau}\,\eta_2)(x)+A\, ({\cal E}^2_{{\tau}}\,\xi_2)(x,\mu)\\
g_3^1(x,0,\mu)&=&A\,  ({\cal E}^3_{\tau}\,\eta_3+Y_3)(x)+A\, ({\cal E}^3_{{\tau}}\,\xi_3 +
{Z}_3)(x,\mu)\\
&\vdots&\\
g_{j}^1(x,0,\mu)&=&A\,  ({\cal E}^j_{\tau}\,\eta_j+Y_j)(x)+A\, ({\cal E}^{j}_{\tau}\,\xi_{j} + {Z}_{j})(x,\mu)\\
&\vdots&
\end{array}
\label{gformhopf}
\end{equation}
In (\ref{gformhopf}), ${\cal E}^j_{\tau}$ are as in (\ref{Ecaldefmulthopf}),
$A$ is the $\mathbb{T}^p$ averaging operator (\ref{Adef}), and 
$Y_j(x)$ and $Z_j(x,\mu)$ are the extra contributions to the terms of order $j$ coming
from the lower order $(<j)$ changes of variables, and $Z_j(x,0)=0$.  Hence, the terms
$Y_j$ and $Z_j$ are
completely determined once the normalizing procedure arrives at order
$j$.  More precisely, $Y_j$ is determined explicitly in terms of
$\eta_2,\ldots,\eta_{j-1},U^i_2,\ldots,U^i_{j-1},\,i=1,2$ and $Z_j$ is determined
explicitly in terms of $\eta_2,\ldots,\eta_{j-1},\xi_2,\ldots,\xi_{j-1},
U^i_2,\ldots,U^i_{j-1}, W^i_2,\ldots,W^i_{j-1},\,i=1,2$. 
Taking into account (\ref{gformhopf}) and using the convention
$Y_2=0$, ${Z}_2=0$, the center manifold equations
  (\ref{y6}) are $\mathbb{T}^p$-equivariant, and in polar coordinates,
  the uncoupled radial part (truncated at order $\ell$) is of the form
\[
\dot{\rho}=\sum_{j=2}^{\ell}\,\left[({\cal N}^j_{\tau}\,\eta_j+(\Pi\circ A)(Y_j))(\rho)+
({\cal
  N}^{j}_{{\tau}}\,\xi_{j}+(\Pi\circ A)({Z}_{j}))(\rho,\mu)\right].
\]
Thus, using Proposition \ref{dsdecomMultipleHopfinv}, we get the desired result if we set 
\begin{equation}
\eta_{j}=\left({\cal
    N}^{j}_{{\tau}}\right)^{-1}(h_{j}-(\Pi\circ
    A)({Y}_{j})),\,\,\,\,\,\,\,
\xi_j=\left({\cal N}^j_{\tau}\right)^{-1}(q_j-(\Pi\circ A)(Z_j)),
\,\,\,j=2,\ldots,\ell.
\label{recipe}
\end{equation}
\hfill\qed

Theorem \ref{mainthm1} has an important interpretation in terms of
the singularity and unfolding theory of scalar delay-differential equations.
Suppose (\ref{y1}) satisfies the hypotheses of Theorem \ref{mainthm1}.
Let $h(\rho)$ be any given (parameter independent) element
of ${\cal H}^{d}_{\ell}(\mathbb{R}^d,\mathbb{Z}_{2,p})$, $\ell\geq 2$.  Then
Theorem \ref{mainthm1} implies that (under generic conditions on
$\tau_1,\ldots,\tau_d$) there exists an unparameterized
nonlinear polynomial delay-differential equation 
\begin{equation}
\dot{z}(t)=L_0\,z_t+\eta(z(t+\tau_1),\ldots,z(t+\tau_d))
\label{unpertpd}
\end{equation}
whose dynamics on a center manifold up to order $\ell$ have
as uncoupled radial equations $\dot{\rho}=h(\rho)$.  
Therefore, generically, {\em any finitely-determined singularity
  within the space of $\mathbb{Z}_{2,p}$-equivariant radial equations
can be realized by an appropriate choice of $\eta$ in (\ref{unpertpd})}.

Now, suppose
that $\tilde{h}\in {\cal H}^{d+s}_{\ell}(\mathbb{R}^d,\mathbb{Z}_{2,p})$ is
an
{\em equivariant unfolding} of the finitely-determined singularity $h$
above, i.e. $\tilde{h}(\rho,\mu)$ is such that
$\tilde{h}(\rho,0)=h(\rho)$.  
Then $q(\rho,\mu)\equiv \tilde{h}(\rho,\mu)-h(\rho)$ is an element of ${\cal P}^{d+s}_{\ell}(\mathbb{R}^d,\mathbb{Z}_{2,p})$.
Theorem \ref{mainthm1} implies that there exists a parameterized
nonlinear polynomial delay-differential equation of the form
\begin{equation}
\dot{z}_t=L_0\,z_t+\eta(z(t+\tau_1),\ldots,z(t+\tau_d))+\xi(z(t+\tau_1),\ldots,z(t+\tau_d),\mu),
\label{pertpd}
\end{equation}
with $\xi(z(t+\tau_1),\ldots,z(t+\tau_d),0)=0$,
whose dynamics on a center manifold up to order $\ell$ have as
uncoupled radial equations $\dot{\rho}=\tilde{h}(\rho,\mu)=h(\rho)+q(\rho,\mu)$. 
Therefore, the unfolding {\em
  $\eta(z(t+\tau_1),\ldots,z(t+\tau_d))+\xi(z(t+\tau_1),\ldots,z(t+\tau_d),\mu)$
of $\eta$
realizes the unfolding $\tilde{h}(\cdot,\mu)$ of the singularity
$h$ on the center manifold}.

In the theory of classification of singularities of equivariant vector fields
\cite{GSSI,GSSII}, one often defines a suitable equivalence relation on
a given space of vector fields (by requiring preservation of certain local
qualitative features of the flow associated to the vector field), and then classifies the
equivalence classes in terms of a (hopefully finite)
set of conditions of the Taylor coefficients of the vector field.
One then wishes to characterize the ``likelihood'' of a given
singularity, $f$,
by computing its {\em codimension}, which roughly speaking, is the
codimension of the equivalence orbit through $f$.
Finally, one
then uses this idea of codimension to construct a versal unfolding of
the singularity (perturbing in transversal directions to the
equivalence orbit).

Suppose
$f:\mathbb{R}^d\longrightarrow\mathbb{R}^d$ is a smooth vector field
vanishing at the origin and equivariant with respect to the group
$\mathbb{Z}_{2,p}$ previously defined.  Typically, the computation of codimension of
the singularity $f$ is done by first identifying a polynomial $h\in
{\cal H}^d_{\ell}(\mathbb{R}^d,\mathbb{Z}_{2,p})$ (for some suitable
$\ell$) which is equivalent to $f$ regardless of the Taylor
coefficients of $f$ of order greater than $\ell$.  Then, one
constructs the tangent space 
within ${\cal H}^d_{\ell}(\mathbb{R}^d,\mathbb{Z}_{2,p})$
to the equivalence orbit
of $h$ through $h$, $T_h\subset {\cal H}^d_{\ell}(\mathbb{R}^d,\mathbb{Z}_{2,p})$.  Finally one finds a complementary subspace
$C_h\subset {\cal
  H}^d_{\ell}(\mathbb{R}^d,\mathbb{Z}_{2,p})$
such that ${\cal
  H}^d_{\ell}(\mathbb{R}^d,\mathbb{Z}_{2,p})=T_h\,\oplus\,C_h$.  The
codimension of $f$ is then identified with the dimension of $C_h$
(i.e. $\mbox{\rm codim}\,f\equiv \mbox{\rm codim}\,T_h=\mbox{\rm dim}\,C_h$).  We
say that $f$ is {\em generic} with respect to the equivalence relation
if the codimension of $f$ is zero.
The following theorem
addresses this issue within the context of realizability.
\begin{thm}
Consider the RFDE (\ref{y1}) in the unparametrized ($s=0$) case, and let $\Lambda_0$ denote the set of
solutions of (\ref{ychar}) with zero real part.  Suppose that
Hypothesis \ref{spectralhyp} is satisfied. 
Suppose that the nonlinear term $F(z_t)$
is of the general form 
\[
F(z_t)=\eta(z(t+\tau_1),\ldots,z(t+\tau_d)),
\]
where $\eta$ is smooth.  Then the local dynamics of (\ref{y1}) near the
origin on an invariant center manifold
can be described by a system of ordinary differential
equations on $\mathbb{R}^{\kappa}$.  Moreover, this ODE system can be
brought into $\mathbb{T}^p$-equivariant normal form to any desired
order $\ell$, and the resulting (truncated at order $\ell$) normal
form can be uncoupled into two sub-systems
\begin{eqnarray}
\dot{\rho}&=&h(\rho\,;\,\eta,\tau)\label{1st2}\\[0.15in]
\dot{\theta}&=&k(\rho\,;\,\eta,\tau),\label{2nd2}
\end{eqnarray}
where $\tau=(\tau_1,\ldots,\tau_d)\in\mathbb{R}^d$,
$h(\cdot\,;\,\eta,\tau)\in {\cal
  H}^d_{\ell}(\mathbb{R}^d,\mathbb{Z}_{2,p})$ and
$k(\cdot\,;\,\eta,\tau):\mathbb{R}^p\longrightarrow\mathbb{R}^p$.
For given $\tau\in\mathbb{R}^d$, consider the following mapping:
\[
\begin{array}{rcl}
{\cal F}_{\tau}:{\cal H}^d_{\ell}(\mathbb{R})&\longrightarrow& {\cal
  H}^{d}_{\ell}(\mathbb{R}^d,\mathbb{Z}_{2,p})\\[0.15in]
\eta&\longmapsto&{\cal F}_{\tau}(\eta)=h(\cdot\,;\,\eta,\tau),
\end{array}
\]
where $h$ is as in (\ref{1st2}).
Then there is an open and dense set ${\cal U}\subset\mathbb{R}^d$,
such that for all $\tau\in {\cal U}$, ${\cal F}_{\tau}$ is a
submersion.
Consequently, if ${\cal M}\subset {\cal
  H}^d_{\ell}(\mathbb{R}^d,\mathbb{Z}_{2,p})$ is a smooth manifold,
then for all $\tau\in {\cal U}$, ${\cal
    F}_{\tau}^{-1}({\cal M})$ is a smooth submanifold of ${\cal
  H}^d_{\ell}(\mathbb{R})$, and
\begin{equation}
\mbox{\rm codim}\,{\cal
    F}_{\tau}^{-1}({\cal M})=
\mbox{\rm codim}\,{\cal M}
\label{codimeq}
\end{equation}
\label{mainthm3}
\end{thm}
\proof
The fact that the center manifold equations are given by (\ref{1st2})
and (\ref{2nd2}) has already been proved.  

The mapping ${\cal F}_{\tau}$ is computable similarly to
(\ref{recipe}): if $\eta=\sum_{j=2}^{\ell}\,\eta_j$, with $\eta_j\in
H^d_{j}(\mathbb{R})$, then
\[
{\cal
  F}_{\tau}(\eta)=h(\cdot\,;\,\eta,\tau)=\sum_{j=2}^{\ell}\,((\Pi\circ
  A\circ {\cal
  E}^j_{\tau})(\eta_j)+(\Pi\circ A)(Y_j)),
\]
where $Y_2=0$ and $Y_j$ is a smooth function of $\eta_2,\ldots,\eta_{j-1}$
for $j>2$.  Thus, if $\zeta=\sum_{j=2}^{\ell}\,\zeta_j$, with $\zeta_j\in
H^d_{j}(\mathbb{R})$, then 
\[
D{\cal F}_{\tau}(\eta)\cdot\zeta=\sum_{j=2}^{\ell}\,\left((\Pi\circ A\circ
{\cal E}^j_{\tau})(\zeta_j)+(\Pi\circ
A)\left(\sum_{i=2}^{\ell}\,Y_{ji}(\eta)\zeta_i\right)\right),
\]
where $Y_{ji}=0$ if $i\geq j$.  From Proposition
  \ref{MultipleHopfinv}, there is is an open and dense set ${\cal
  U}\subset\mathbb{R}^d$ such that for all $\tau\in{\cal U}$,
$D{\cal F}_{\tau}(\eta)$ is onto ${\cal
  H}^d_{\ell}(\mathbb{R}^d,\mathbb{Z}_{2,p})$, and consequently ${\cal
  F}_{\tau}$ is a submersion.  Equation (\ref{codimeq}) follows from
  the transversal mapping theorem \cite{AMR}.
\hfill\qed

The next result states that the number of delays, $d$, shown above to
be sufficient to
realize any arbitrary element of ${\cal
  H}^{d}_{\ell}(\mathbb{R}^d,\mathbb{Z}_{2,p})$, is optimal for that purpose.
\begin{thm}
Consider the RFDE (\ref{y1}) in the unparameterized ($s=0$) case, and let $\Lambda_0$ denote the set of
solutions of (\ref{ychar}) with zero real part.  Suppose that
Hypothesis \ref{spectralhyp} is satisfied. 
Suppose that the nonlinear term $F(z_t)$
is of the general form 
\[
F(z_t)=\eta(z(t+\tau_1),\ldots,z(t+\tau_{d-1})),
\]
where $\eta$ is smooth.  Then the local dynamics of (\ref{y1}) near the
origin on an invariant center manifold
can be described by a system of ordinary differential
equations on $\mathbb{R}^{\kappa}$.  Moreover, this ODE system can be
brought into $\mathbb{T}^p$-equivariant normal form to any desired
order $\ell$, and the resulting (truncated at order $\ell$) normal
form can be uncoupled into two sub-systems
\begin{eqnarray}
\dot{\rho}&=&h(\rho\,;\,\eta,\tau)\label{1st3}\\[0.15in]
\dot{\theta}&=&k(\rho\,;\,\eta,\tau),\label{2nd3}
\end{eqnarray}
where $\tau=(\tau_1,\ldots,\tau_{d-1})\in\mathbb{R}^{d-1}$,
$h(\cdot\,;\,\eta,\tau)\in {\cal
  H}^d_{\ell}(\mathbb{R}^d,\mathbb{Z}_{2,p})$ and
$k(\cdot\,;\,\eta,\tau):\mathbb{R}^p\longrightarrow\mathbb{R}^p$.
For given $\tau\in\mathbb{R}^{d-1}$, consider the following mapping:
\begin{equation}
\begin{array}{rcl}
{\cal F}_{\tau}:{\cal H}^{d-1}_{\ell}(\mathbb{R})&\longrightarrow& {\cal
  H}^{d}_{\ell}(\mathbb{R}^d,\mathbb{Z}_{2,p})\\[0.15in]
\eta&\longmapsto&{\cal F}_{\tau}(\eta)=h(\cdot\,;\,\eta,\tau),
\end{array}
\label{Fdm1}
\end{equation}
where $h$ is as in (\ref{1st3}).  Then there is an integer $\ell_0\geq
2$ such that ${\cal F}_{\tau}$ is not surjective if $\ell\geq\ell_0$.
\label{restrictionsthm}
\end{thm}
\proof
It will be sufficient to show that for fixed $d$, 
\begin{equation}
\frac{\mbox{\rm dim}\,{\cal H}^{d-1}_{\ell}(\mathbb{R})}{\mbox{\rm
    dim}\,{\cal
    H}^d_{\ell}(\mathbb{R}^d,\mathbb{Z}_{2,p})}=O(\ell^{-1})
\label{bigOhell}
\end{equation}
as $\ell\rightarrow\infty$.
First, note that it is well-known \cite{GJ} that for given integers $m\geq 1$
and $\ell\geq 1$, the number of solutions in
non-negative integers for the equation
\[
k_1+\cdots +k_{m}=\ell
\]
is
\[
\left(\begin{array}{c}
m+\ell-1\\m-1\end{array}\right).
\]
Thus,
\begin{equation}
\begin{array}{lll}
\mbox{\rm dim}\,{\cal
  H}^{d-1}_{\ell}(\mathbb{R})&=&{\dps\sum_{j=2}^{\ell}\,\mbox{\rm
  dim}\,H^{d-1}_{j}(\mathbb{R})=
\sum_{j=2}^{\ell}\,\left(\begin{array}{c}
d+j-2\\d-2\end{array}\right)}\\\\
&=&{\dps
  \left(\begin{array}{c}d-1+\ell\\d-1\end{array}\right)-d=O(\ell^{d-1})\,\,\,\,\mbox{\rm
  as}\,\,\,\ell\rightarrow\infty}.
\end{array}
\label{HRdim}
\end{equation}
Using a similar (but slightly lengthier) computation, we can show that
\[
\mbox{\rm dim}\,{\cal
  H}^d_{\ell}(\mathbb{R}^d,\mathbb{Z}_{2,p})=O(\ell^d)\,\,\,\,\mbox{\rm as}\,\,\,\ell\rightarrow\infty,
\]
which establishes (\ref{bigOhell}) and concludes the proof of this theorem.
\hfill\qed

Theorem \ref{restrictionsthm} is important in the problem of
establishing whether or not there are restrictions on the possible phase portraits
for an unfolding of a given singularity $h\in {\cal
  H}^d_{\ell}(\mathbb{R}^d,\mathbb{Z}_{2,p})$ when such an unfolding
arises from center manifold reduction and phase/amplitude decoupling
of a nonlinear delay-differential equation (\ref{y1}).  This theorem
allows one to conclude that, at least for $\ell$ large enough, such
restrictions are likely to occur if the number of delays in the
nonlinearity $\tilde{F}$ in (\ref{y1p}) is less than $d$.  For example,
this question was addressed in \cite{BB} in the context of the
non-resonant double Hopf bifurcation.  In this case, they show that if
the nonlinear part of (\ref{y1}) contains 2 delays, then generically 
any cubic order radial equation (\ref{1st2}) can be realized by
appropriate choice of the nonlinear coefficients in (\ref{y1}) (note
that our Theorem \ref{mainthm1} recovers and generalizes that result).
However, in \cite{BB}, it is also shown
that if the
nonlinear part of (\ref{y1}) depends on only one delay, then not all
equivalence classes of phase portraits in the versal unfolding of the
radial equations (\ref{1st3}) can be attained by variation of the
nonlinear coefficients in (\ref{y1}), for fixed values of $\omega_1$,
$\omega_2$ and $\tau$.  The next example treats this specific case in
the context of Theorem \ref{restrictionsthm}.
\begin{examp}
{\rm
In the case where $\Lambda_0=\{\pm\,i\omega_1,\pm\,i\omega_2\}$, we
get from (\ref{HRdim}) that 
\[
\mbox{\rm dim}\,{\cal H}^1_{\ell}(\mathbb{R})=\ell-1.
\]
It is also easy to show that
if $\ell=2L+j$, where $L\geq 0$ is an integer and $j\in\{0,1\}$, then
\[
\mbox{\rm dim}\,{\cal
  H}^2_{\ell}(\mathbb{R}^2,\mathbb{Z}_{2,2})=L(L+3).
\]
Thus, from $\ell=3$ onward, we have $\mbox{\rm dim}\,{\cal
  H}^1_{\ell}(\mathbb{R})<\mbox{\rm dim}\,{\cal
  H}^2_{\ell}(\mathbb{R}^2,\mathbb{Z}_{2,2})$.
In particular, $\mbox{\rm dim}\,{\cal
  H}^1_{3}(\mathbb{R})=2$ and $\mbox{\rm dim}\,{\cal
  H}^2_{3}(\mathbb{R}^2,\mathbb{Z}_{2,2})=4$.
Therefore, at cubic order, the mapping ${\cal F}_{\tau}$ in (\ref{Fdm1})
\[
\begin{array}{rcl}
{\cal F}_{\tau}:{\cal H}^{1}_{3}(\mathbb{R})&\longrightarrow& {\cal
  H}^{2}_{3}(\mathbb{R}^2,\mathbb{Z}_{2,2})\\[0.15in]
\eta&\longmapsto&{\cal F}_{\tau}(\eta)=h(\cdot\,;\,\eta,\tau),
\end{array}
\]
is not surjective, and so there are elements of ${\cal
  H}^2_{3}(\mathbb{R}^2,\mathbb{Z}_{2,2})$ which can not be realized
  by any element of ${\cal H}^1_3(\mathbb{R})$.  In fact, ${\cal
  F}_{\tau}\left({\cal H}^1_3(\mathbb{R})\right)$ is a two-dimensional
  smooth surface in the 4-dimensional space ${\cal
  H}^2_3(\mathbb{R}^2,\mathbb{Z}_{2,2})$, such that ${\cal
  F}_{\tau}(0)=0$.  Specifically, if we write the general element of
  ${\cal H}^1_3(\mathbb{R})$ as
\[
b_2v^2+b_3v^3
\]
and the general element of ${\cal
  H}^2_3(\mathbb{R}^2,\mathbb{Z}_{2,2})$ as
\[
\begin{array}{c}
(a_{11}\rho_1^2+a_{12}\rho_2^2)\rho_1\\\\
(a_{21}\rho_1^2+a_{22}\rho_2^2)\rho_2,
\end{array}
\]
then the mapping ${\cal F}_{\tau}$ can be represented by the following
mapping from $\mathbb{R}^2$ into $\mathbb{R}^4$:
\begin{equation}
a_{ij}(b_2,b_3)=\alpha_{ij}\,b_2^2+\beta_{ij}\,b_3,\,\,\,\,\,i,j=1,2
\label{Fcoords}
\end{equation}
where the real coefficients $\alpha_{ij}$ and $\beta_{ij}$ are
determined from $\tau$, $\omega_1$ and $\omega_2$.

Note however that the problem of determining whether or not there are
  restrictions is somewhat more subtle than one of surjectivity, since the topological types of
  the possible phase diagrams in the unfolding space for the double
  Hopf bifurcation are determined by the {\em sign} of the cubic
  coefficients in the radial equations (and not their actual values).
In the $(b_2,b_3)$ plane, the zero level sets of the $a_{ij}$ in
  (\ref{Fcoords}) are (at most) four distinct curves (parabolae generically) which
  intersect only at the origin.  Consequently, there are at most four
  distinct open regions in the $(b_2,b_3)$ plane in which the signs of
  the coefficients $a_{ij}$ are constant and non-zero.  It is then
  easy to see that it is impossible to realize the twelve possible
  sign combinations (see \cite{GH}) which characterize the complete
  unfolding space of the double Hopf bifurcation, and so there will be
  restrictions on the phase portraits when the nonlinear terms in
  (\ref{y1}) contain only one delay.
}
\end{examp}

For general RFDEs (i.e. not necessarily delay-differential equations),
we have the following result on realization of unfoldings:
\begin{thm}
Consider the general nonlinear RFDE
\begin{equation}
\dot{z}(t)=L_0\,z_t+N(z_t)
\label{genunpertrfde}
\end{equation}
where $L_0:C\rightarrow\mathbb{R}$ is
a bounded linear
operator from $C\equiv C\left(  \left[  -r,0\right]  ,\mathbb{R}\right)$ into $\mathbb{R}$, and $N$ is
a smooth function from $C$ into $\mathbb{R}$,
with $N(0)=0$, 
$DN(0)=0$.
Let $\Lambda_0$ denote the set of
solutions of (\ref{ychar}) with zero real part and suppose that
Hypothesis \ref{spectralhyp} is satisfied. 
Then the local dynamics of (\ref{genunpertrfde})
near the origin on an invariant
center manifold can be described by a system of ordinary differential
equations on $\mathbb{R}^{\kappa}$.  Moreover, this ODE system can be
brought into $\mathbb{T}^p$-equivariant normal form to any desired
order $\ell$, and the resulting (truncated at order $\ell$) normal
form can be uncoupled into an uncoupled $d$-dimensional system and a
$p$-dimensional system
\begin{eqnarray}
\dot{\rho}&=&h(\rho\,;\,N)\label{1st}\\[0.15in]
\dot{\theta}&=&k(\rho\,;\,N),\label{2nd}
\end{eqnarray}
where for given $N$, $h(\cdot\,;\,N)$ is some element of ${\cal
  H}^{d}_{\ell}(\mathbb{R}^d,\mathbb{Z}_{2,p})$, and
  $k(\cdot\,;\,N):\mathbb{R}^p\longrightarrow\mathbb{R}^p$.
Let $\tilde{h}(\rho,\mu)$ be an $s$-parameter equivariant unfolding of $h$ of degree at
  most $\ell$, i.e.
$\tilde{h}\in {\cal H}^{d+s}_{\ell}(\mathbb{R}^d,\mathbb{Z}_{2,p})$ and
  $\tilde{h}(\cdot,0)=h(\cdot\,;\,N)$.
Then there exists an $s$-parameter unfolding of (\ref{genunpertrfde})
  of the form
\begin{equation}
\dot{z}(t)=L_0(z_t)+N(z_t)+\xi(z(t+\tau_1),\ldots,z(t+\tau_d),\mu)
\label{genunfrfde}
\end{equation}
(where $\tau=(\tau_1,\ldots,\tau_d)\in\mathbb{R}^d$, and
$\xi\in\widehat{{\cal W}}^{d+s}_{\ell}(\mathbb{R})$ vanishes at $\mu=0$)
which realizes the unfolded radial equations
\[
\dot{\rho}=\tilde{h}(\rho,\mu)
\]
on an invariant center manifold for (\ref{genunfrfde}).
\label{mainthm2}
\end{thm}
\proof 
Choosing $W^1_j=0$, $W^2_j=0$ and choosing $U^1_j$ and $U^2_j$
appropriately in (\ref{dsdecomnfcv}), 
the center manifold equations for (\ref{genunpertrfde}) truncated at
order $\ell$ are equivalent to
(\ref{1st}) and (\ref{2nd}).  

Using 
Theorem \ref{thmdsdecomnf},
for arbitrary $\xi=\sum_{j=2}^{\ell}\,\xi_j\in\widehat{{\cal W}}^{d+s}_{\ell}(\mathbb{R})$,
 there is a sequence of near-identity changes of variables
 (\ref{dsdecomnfcv}) (with $U^1_j$ and $U^2_j$ as above) for which the
 uncoupled radial part of the center manifold equations for
 (\ref{genunfrfde}) truncated at
 order $\ell$ are
\[
\dot{\rho}=h(\rho\,;\,N)+\sum_{j=2}^{\ell}\,({\cal
 N}^j_{\tau}\,\xi_j+(\Pi\circ A)(Z_j))
\]
where $(\Pi\circ A)(Z_j)$ is some known element of
${P}^{d+s}_{j}(\mathbb{R}^d,\mathbb{Z}_{2,p})$.
The conclusion follows from setting
\[
\xi_j=\left({\cal N}^j_{\tau}\right)^{-1}(\tilde{h}-h-(\Pi\circ
A)(Z_j)).
\]
\hfill
\qed

\Section{The \mbox{\boldmath $\pm\,i\omega$}, \mbox{\boldmath $(\pm\, i\omega_1,\pm\,i\omega_2)$} and
  \mbox{\boldmath $(0,\pm\,i\omega)$} Singularities}

Our results in Theorems \ref{mainthm1}, \ref{mainthm3} and \ref{mainthm2}
allow us to recover some previous results on realizability and
(lack of) restrictions for Hopf bifurcation, non-resonant double Hopf
bifurcation, and the $(0,\pm i\omega)$ singularity in
scalar RFDEs \cite{FM96,BB}.
\begin{cor}[Theorem 1 of \cite{FM96}]
Consider the RFDE (\ref{y1}) in the unparameterized case
\begin{equation}
\dot{z}(t)=L_0z_t+F(z_t),
\label{y30}
\end{equation}
such that the
characteristic equation (\ref{ychar}) has simple purely imaginary
roots $\pm i\omega\neq 0$ and no other
roots on the imaginary axis (simple Hopf bifurcation).  If
\begin{equation}
F(z_t)=A_2(z(t+\tau))^2+A_3(z(t+\tau))^3,\,\,\,\tau\in [-r,0]
\label{FhopfFM}
\end{equation}
then the
uncoupled radial part of the center manifold equations to cubic order are
\begin{equation}
\dot{\rho}=a\rho^3,
\label{hopfnfFM}
\end{equation}
where $a=a(A_2,A_3;\tau,\omega)$.  Generically, the non-degeneracy condition
$a\neq 0$ is satisfied.
In fact, for any $a\in\mathbb{R}$, (\ref{hopfnfFM}) can be realized
with $A_2=0$ for an appropriate choice of $A_3$ in (\ref{FhopfFM}).
Furthermore, in the case $a\neq 0$, the versal unfolding
\[
\dot{\rho}=\mu\rho+a\rho^3
\]
of (\ref{hopfnfFM}) is generically realized (modulo a rescaling of the
parameter) by the following unfolding
of (\ref{y30})
\[
\dot{z}(t)=L_0z_t+F(z_t)+\mu\,z(t+\tau).
\]
\end{cor}
\begin{cor}[Theorem 3.1(1) of \cite{BB}]
Consider the RFDE (\ref{y1}) in the unparameterized case
\begin{equation}
\dot{z}(t)=L_0z_t+F(z_t),
\label{y40}
\end{equation}
such that the
characteristic equation (\ref{ychar}) has simple non-resonant purely imaginary
roots $\pm i\omega_1, \pm i\omega_2$, and no other
roots on the imaginary axis (non-resonant double Hopf bifurcation).  If
\begin{equation}
\begin{array}{lll}
F(z_t)&=&A_{20}(z(t+\tau_1))^2+A_{11}z(t+\tau_1)z(t+\tau_2)+A_{02}(z(t+\tau_2))^2+\\[0.15in]
&&A_{30}(z(t+\tau_1))^3+A_{21}(z(t+\tau_1))^2z(t+\tau_2)+A_{12}z(t+\tau_1)(z(t+\tau_2))^2+\\[0.15in]
&&A_{03}(z(t+\tau_2))^3,
\end{array}
\label{F2hopfBB}
\end{equation}
where $\tau_1,\tau_2\in [-r,0]$,
then the
uncoupled radial part of the center manifold equations to cubic order are
\begin{equation}
\begin{array}{lll}
\dot{\rho_1}&=&(a_{11}\rho_1^2+a_{12}\rho_2^2)\rho_1\\[0.15in]
\dot{\rho_2}&=&(a_{21}\rho_1^2+a_{22}\rho_2^2)\rho_2,
\label{2hopfnfBB}
\end{array}
\end{equation}
where
$a_{ij}=a_{ij}(A_{20},A_{11},A_{02},A_{30},A_{21},A_{12},A_{03};\tau_1,\tau_2,\omega_1,\omega_2)$.
Generically, the non-degeneracy condition
$a_{11}a_{22}-a_{21}a_{12}\neq 0$ is satisfied.
In fact, for any $a_{11}, a_{12}, a_{21}, a_{22}\in\mathbb{R}$, (\ref{2hopfnfBB}) can be realized
with $A_{20}=A_{11}=A_{02}=0$ for an appropriate choice of $A_{30},
A_{21}, A_{12}, A_{03}$ in (\ref{F2hopfBB}).
Furthermore, in the case $a_{11}a_{22}-a_{12}a_{21}\neq 0$, the versal unfolding
\begin{equation}
\begin{array}{lll}
\dot{\rho_1}&=&(\mu_1+a_{11}\rho_1^2+a_{12}\rho_2^2)\rho_1\\[0.15in]
\dot{\rho_2}&=&(\mu_2+a_{21}\rho_1^2+a_{22}\rho_2^2)\rho_2
\end{array}
\label{2hopfnfBBunf}
\end{equation}
of (\ref{2hopfnfBB}) is generically realized (modulo a linear change of
parameters) by the following unfolding
of (\ref{y40})
\[
\dot{z}(t)=L_0z_t+F(z_t)+\mu_1\,z(t+\tau_1)+\mu_2\,z(t+\tau_2).
\]
\label{corBB}
\end{cor}
\begin{rmk}
We would like to clarify the statement ``modulo a linear change of
parameters'' in the preceding Corollary.  According to the notation we
have established in this paper, we have
\[
V_2^{2+2}(\mathbb{R})=\widehat{V}_2^{2+2}(\mathbb{R})=W_2^{2+2}(\mathbb{R})=\widehat{W}_2^{2+2}(\mathbb{R})=\mbox{\rm
  span}\,\{\,\mu_1\,v_1,\,\mu_1\,v_2,\,\mu_2\,v_1,\,\mu_2\,v_2\,\}
\]
and
\[
H_2^{2+2}(\mathbb{R}^2,\mathbb{Z}_{2,2})=P_2^{2+2}(\mathbb{R}^2,\mathbb{Z}_{2,2})=\mbox{\rm
  span}\,
\left\{\,\mu_1\,\left(\begin{array}{c}\rho_1\\0\end{array}\right),\,\mu_1\,\left(\begin{array}{c}0\\\rho_2\end{array}\right),\,
\mu_2\,\left(\begin{array}{c}\rho_1\\0\end{array}\right),\,\mu_2\,\left(\begin{array}{c}0\\\rho_2\end{array}\right)\,\right\}.
\]
From Proposition \ref{MultipleHopfinv}, the mapping
\[
{\cal N}^2_{\tau} : V_2^{2+2}(\mathbb{R})\longrightarrow
H_2^{2+2}(\mathbb{R}^2,\mathbb{Z}_{2,2})
\]
is generically invertible.  Since the mapping ${\cal
  N}^2_{\tau}$ does not have any effect on the parameters $\mu_1$ and $\mu_2$, generically we have
\[
\begin{array}{rcl}
({\cal
  N}^2_{\tau})^{-1}\left(\mu_j\,\left(\begin{array}{c}\rho_1\\0\end{array}\right)\right)&=&\mu_j\,(m_{11}v_1+m_{12}v_2),\,\,\,j=1,2\\
&&\\
({\cal
  N}^2_{\tau})^{-1}\left(\mu_j\,\left(\begin{array}{c}0\\\rho_2\end{array}\right)\right)&=&\mu_j\,(m_{21}v_1+m_{22}v_2),\,\,\,j=1,2,\end{array}
\]
and it follows that the $2\times 2$ matrix $M=(m_{ij})$ is
invertible.  Consequently, the unfolding (\ref{2hopfnfBBunf}) of
(\ref{2hopfnfBB}) is realized by the following unfolding of
(\ref{y40})
\[
\dot{z}(t)=L_0z_t+F(z_t)+\mu_1\,(m_{11}z(t+\tau_1)+m_{12}z(t+\tau_2))+\mu_2\,(m_{21}z(t+\tau_1)+m_{22}z(t+\tau_2)).
\]
We get the conclusion of Corollary \ref{corBB} by performing the
linear change of parameters
\[
\tilde{\mu}_1=m_{11}\mu_1+m_{21}\mu_2,\,\,\,\,\,\,\,\,
\tilde{\mu}_2=m_{12}\mu_1+m_{22}\mu_2
\]
and dropping the tildes.
\end{rmk}
\begin{cor}[Theorem 2 of \cite{FM96}]
Consider the RFDE (\ref{y1}) in the unparameterized case
\begin{equation}
\dot{z}(t)=L_0z_t+F(z_t),
\label{y50}
\end{equation}
such that the
characteristic equation (\ref{ychar}) has simple purely imaginary
roots $\pm i\omega\neq 0$, a simple root at $0$, and no other
roots on the imaginary axis (interaction of a simple bifurcation and a
Hopf bifurcation).
If
\begin{equation}
F(z_t)=A_{20}(z(t+\tau_1))^2+A_{11}z(t+\tau_1)z(t+\tau_2)+A_{02}(z(t+\tau_2))^2
\label{FsshopfFM}
\end{equation}
where $\tau_1,\tau_2\in [-r,0]$,
then the
uncoupled radial part of the center manifold equations to quadratic order are
\begin{equation}
\begin{array}{lll}
\dot{\rho_0}&=&b_{1}\rho_0^2+b_{2}\rho_1^2\\[0.15in]
\dot{\rho_1}&=&a_1\rho_0\rho_1,
\label{sshopfnfFM}
\end{array}
\end{equation}
where the coefficients $a_1$, $b_1$ and $b_2$ are functions of
$(A_{20},A_{11},A_{02};\tau_1,\tau_2,\omega)$.  Generically, the
non-degeneracy conditions $a_1\neq 0$, $b_1\neq 0$, $b_2\neq 0$ and
$a_1\neq b_2$ are satisfied.  
\end{cor}
\begin{cor}
Consider the singularity (\ref{sshopfnfFM}) in the non-degenerate case
$a_1\neq 0$, $b_1\neq 0$, $b_2\neq 0$ and
$a_1\neq b_2$.  Then the following {\em Langford unfolding}
\cite{L79} of (\ref{sshopfnfFM}) in the transcritical case
\begin{equation}
\begin{array}{lll}
\dot{\rho_0}&=&\mu_1\rho_0+b_{1}\rho_0^2+b_{2}\rho_1^2\\[0.15in]
\dot{\rho_1}&=&\mu_2\rho_1+a_1\rho_0\rho_1,
\label{sshopfnfFMunf}
\end{array}
\end{equation}
is generically realized (modulo a linear change of parameters) by the
following unfolding of (\ref{y50})
\[
\dot{z}(t)=L_0z_t+F(z_t)+\mu_1\,z(t+\tau_1)+\mu_2\,z(t+\tau_2).
\]
\end{cor}

\Section{Conclusions}
We have established a framework for the realizability problem for
scalar RFDEs
which exploits fully the toroidal equivariance of normal forms of bifurcations
associated with purely imaginary eigenvalues.  This has allowed us to
recover and significantly generalize recent results of Faria and
Magalh$\tilde{\mbox{\rm a}}$es \cite{FM96} and of Buono and B\'elair
\cite{BB}.  As mentioned in the Introduction, it is important for
modelers using RFDEs to be able to accurately assess the range of
possible dynamics accessible within their models.  For this purpose,
this paper gives a thorough analysis of this question in the case
where the model is a nonlinear delay-differential equation undergoing 
non-resonant multiple Hopf bifurcation or transcritical/non-resonant
multiple Hopf interaction.  Specifically, we split the dynamics of the
normal form into components which are normal to the orbits of a torus
group, and components which are tangent to these group orbits.  Sharp
estimates on the number of delays are then given for the
realizability of the normal ``radial'' part of the normal form by nonlinear
delay-differential equations.
The case of saddle-node/non-resonant
multiple Hopf interaction will be treated using similar techniques in a subsequent paper \cite{CL}.

The generalizations we have achieved in our paper are twofold.
First, we can treat within a unified framework the general case of $p$
non-resonant Hopf eigenvalues and the interaction between simple steady-state
bifurcation and $p$ non-resonant Hopf bifurcation.
Second, in contrast to \cite{FM96} and \cite{BB} where only the generic
(non-degenerate) cases are treated, we can treat the general
finitely-determined case (whether degenerate or not) and its
unfoldings, also within a unified framework.  Note that
in parameterized families of vector fields with sufficiently many parameters, it becomes possible to violate any
specified non-degeneracy condition which is expressed in terms of the Taylor
coefficients of the vector field up to some finite order.
Therefore, it becomes desirable to have a framework in which these
degenerate cases and their unfoldings can be systematically treated.
Our results provide such a framework.  

Open problems of interest related to this analysis and worthy of
further investigation are
\begin{itemize}
\item relaxing the restriction to scalar RFDEs in order to consider
  $n>1$ dimensional systems of RFDEs
\item incorporating resonances in the purely imaginary eigenvalues and
  repeated eigenvalues with Jordan blocks.
\end{itemize}

\vspace*{0.25in}
\noindent
{\Large\bf Acknowledgments}

\vspace*{0.2in}
This research is partly supported by the
Natural Sciences and Engineering Research Council of Canada in the
form of a Discovery Grant, and by a Premier's Research Excellence
Award from the Ontario Ministry of Economic Development and Trade and the University of Ottawa.

\appendix

\Section{Proof of Proposition~\ref{prop_enf2}}

Let $f$ be a given element of
$H^{\kappa+s}_{\ell}(\mathbb{R}^{\kappa})$, and consider
$\tilde{f}=(f,0)\in H^{\kappa+s}_{\ell}(\mathbb{R}^{\kappa+s})$.  From
Proposition \ref{prop_enf1}, there exists $\tilde{h}=(h_1,h_2)\in
H^{\kappa+s}_{\ell}(\mathbb{R}^{\kappa+s})$ and a unique
$\tilde{g}=(g_1,g_2)\in
H^{\kappa+s}_{\ell}(\mathbb{R}^{\kappa+s},\Gamma)$ such that
\begin{equation}
(f,0)={\cal L}_{\tilde{B}}(h_1,h_2)+(g_1,g_2).
\label{hom_extend}
\end{equation}
Now, ${\cal L}_{\tilde{B}}(h_1,h_2)=({\cal L}_{B}h_1,D_xh_2Bx)$, so it
follows that $g_2=-D_xh_2Bx$.  Consequently, (\ref{hom_extend}) can be
rewritten as
\[
(f,0)=({\cal L}_{B}h_1,0)+(g_1,0),
\]
and thus
\begin{equation}
f={\cal L}_{B}h_1+g_1,
\label{fdecom}
\end{equation}
where $g_1\in H^{\kappa+s}_{\ell}(\mathbb{R}^{\kappa},\mathbb{T}^p)$.
So
$H^{\kappa+s}_{\ell}(\mathbb{R}^{\kappa})=H^{\kappa+s}_{\ell}(\mathbb{R}^{\kappa},\mathbb{T}^p)+\mbox{\rm
  range}
\,{\cal L}_{B}$.  Suppose $f=0$ in (\ref{fdecom}), then it is easy
to see that ${\cal L}_{\tilde{B}}(h_1,0)+(g_1,0)=(0,0)$, and from
Proposition \ref{prop_enf1}, it follows that $g_1=0$ and ${\cal
  L}_{B}h_1=0$.  Therefore,
\[
H^{\kappa+s}_{\ell}(\mathbb{R}^{\kappa})=H^{\kappa+s}_{\ell}(\mathbb{R}^{\kappa},\mathbb{T}^p)\oplus\mbox{\rm
  range}
\,{\cal L}_{B}
\]
\hfill
\qed

\Section{Proof of Proposition~\ref{prop_Adecomposition}}

For a given $f\in H^{\kappa+s}_{\ell}(\mathbb{R}^{\kappa})$, let
$g=Af$; then
\[
\begin{array}{lll}
(Ag)(x,\mu)&=&
{\displaystyle\int_{\Gamma_0}\,\tilde{\gamma}\,g(\tilde{\gamma}^{-1}x,\mu)\,d\tilde{\gamma}=\int_{\Gamma_0}\,\tilde{\gamma}\,\left(\int_{\Gamma_0}\,\gamma\,f(\gamma^{-1}\tilde{\gamma}^{-1}x,\mu)\,d\gamma\right)\,d\tilde{\gamma}}\\[0.15in]
&=&{\displaystyle\int_{\Gamma_0}\,\left(\int_{\Gamma_0}\,\tilde{\gamma}\gamma\,f((\tilde{\gamma}\gamma)^{-1}x,\mu)\,d\gamma\right)\,d\tilde{\gamma}}\\[0.15in]
&=&{\displaystyle\int_{\Gamma_0}\,\left(\int_{\Gamma_0}\,\gamma\,f(\gamma^{-1}x,\mu)\,d\gamma\right)\,d\tilde{\gamma}=\int_{\Gamma_0}\,\gamma\,f(\gamma^{-1}x,\mu)\,d\gamma}\\[0.15in]
&=&(Af)(x,\mu),
\end{array}
\]
where the second to last line holds because of the translation
invariance and the normalization of the Haar integral.
So $A$ is a projection.

Now, let $f\in\mbox{\rm range}\,A$, then $Af=f$, i.e.
\[
f(x,\mu)=\int_{\Gamma_0}\,\gamma\,f(\gamma^{-1}x,\mu)\,d\gamma.
\]
So, for any $\sigma\in\Gamma_0$, we have
\[
\begin{array}{lll}
\sigma\,f(\sigma^{-1}x,\mu)&=&
{\displaystyle\sigma\,\int_{\Gamma_0}\,\gamma\,f(\gamma^{-1}\sigma^{-1}x,\mu)\,d\gamma
  =
  \int_{\Gamma_0}\,\sigma\gamma\,f((\sigma\gamma)^{-1}x,\mu)\,d\gamma}\\[0.15in]
&=&{\displaystyle\int_{\Gamma_0}\,\gamma\,f(\gamma^{-1}x,\mu)\,d\gamma =
  f(x)}.
\end{array}
\]
Therefore,  $f\in
H^{\kappa+s}_{\ell}(\mathbb{R}^{\kappa},\mathbb{T}^p)$.  On the other
hand, if  
$f\in
H^{\kappa+s}_{\ell}(\mathbb{R}^{\kappa},\mathbb{T}^p)$, then
\[
(Af)(x,\mu)=\int_{\Gamma_0}\,\gamma\,f(\gamma^{-1}x,\mu)\,d\gamma =
\int_{\Gamma_0}\,f(x,\mu)\,d\gamma = f(x,\mu),
\]
so $f\in\mbox{\rm range}\,\, A$.  
This establishes (\ref{Arange}).
We now establish (\ref{Aker}).  Since $A$ is a projection, then
\[
H^{\kappa+s}_{\ell}(\mathbb{R}^{\kappa})=\mbox{\rm
  range}\,A\oplus\mbox{\rm ker}\,A.
\]
From Proposition \ref{prop_enf2}, we conclude that $\mbox{\rm dim
  ker}\,A=\mbox{\rm dim range}\,{\cal L}_{B}$.  Thus, we need only show that
$\mbox{\rm range}\,{\cal L}_{B}\subset\mbox{\rm ker}\,A$.  In order to show this, we will
  need the following
\begin{lemma}
Let $g:\Gamma_0\longrightarrow\mathbb{R}^{\kappa}$ be a continuous
function, then
\[
\int_{\Gamma_0}\,g(\gamma)\,d\gamma =
\lim_{T\rightarrow\infty}\,\frac{1}{T}\,\int_0^T\,g(e^{Bs})\,ds.
\]
\label{averaging_lemma}
\end{lemma}
\proof
For a given $q\in\{1,2,\ldots,p\}$, consider the rotation matrix
\[
R_q(\theta)=\mbox{\rm diag}(e^{i\omega_q\theta},e^{-i\omega_q\theta})
\]
which is $T_q\equiv\,2\pi/\omega_q$-periodic in $\theta$.
Then
$\tilde{R}_q(\theta)\equiv R_q(\theta T_q/(2\pi))$ is $2\pi$-periodic
in $\theta$.  By hypothesis, the set
\[
\left\{\frac{2\pi}{T_1},\ldots,\frac{2\pi}{T_p}\right\}
\]
is algebraically independent.  Let $\mbox{\bf T}^p\equiv [0,2\pi]^p$, then we can
parameterize $\Gamma_0$ as follows:
\[
h:\mbox{\bf T}^p\longrightarrow\Gamma_0
\]
\[
h(\theta_1,\ldots,\theta_p)=\left\{\begin{array}{lcl}
\mbox{\rm
  diag}(\tilde{R}_1(\theta_1),\ldots,\tilde{R}_p(\theta_p))\,\,\,&\mbox{\rm if}&\kappa=2p\\[0.15in]
\mbox{\rm
  diag}(1,\tilde{R}_1(\theta_1),\ldots,\tilde{R}_p(\theta_p))&\mbox{\rm if}&\kappa=2p+1.\end{array}\right.
\]
Define $\tilde{g}:\mbox{\bf T}^p\longrightarrow\mathbb{R}^{\kappa}$ by
$\tilde{g}(\theta_1,\ldots,\theta_p)=(g\circ
h)(\theta_1,\ldots,\theta_p)$.  Obviously, $\tilde{g}$ is
$2\pi$-periodic in each of its entries, and 
\[
\frac{1}{(2\pi)^p}\,\int_0^{2\pi}\cdots
\int_0^{2\pi}\,\tilde{g}(\theta_1,\ldots,\theta_p)\,d\theta_1\cdots
d\theta_p=\int_{\Gamma_0}\,g(\gamma)\,d\gamma.
\]
Noting that $h\left(\frac{2\pi s}{T_1},\ldots,\frac{2\pi
      s}{T_p}\right)=e^{Bs}$ and using Lemma 4.1, P.430 of \cite{CH},
      we get that
\[
\begin{array}{lll}
{\displaystyle\frac{1}{(2\pi)^p}\,\int_0^{2\pi}\cdots
\int_0^{2\pi}\,\tilde{g}(\theta_1,\ldots,\theta_p)\,d\theta_1\cdots
d\theta_p}&=&{\displaystyle
\lim_{T\rightarrow\infty}\,\frac{1}{T}\,\int_0^T\,\tilde{g}\left(\frac{2\pi
    s}{T_1},\ldots,\frac{2\pi s}{T_p}\right)\,ds}\\[0.15in]
&=&{\displaystyle\lim_{T\rightarrow\infty}\,\frac{1}{T}\,\int_0^T\,g(e^{Bs})\,ds},
\end{array}
\]
which yields the desired result.
\hfill\qed

Now, let $f\in\mbox{\rm range}\,{\cal L}_{B}$; then there exists $g\in
H^{\kappa+s}_{\ell}(\mathbb{R}^{\kappa})$ such that
\[
D_xg(x,\mu)Bx-Bg(x,\mu)=f(x,\mu),\,\,\,\forall\,(x,\mu)\in\mathbb{R}^{\kappa+s}.
\]
Therefore, using Lemma \ref{averaging_lemma}, we get
\[
\begin{array}{l}
{\displaystyle
(Af)(x,\mu)=\int_{\Gamma_0}\,\gamma\,f(\gamma^{-1}x,\mu)\,d\gamma =
\lim_{T\rightarrow\infty}\,\frac{1}{T}\int_0^T\,e^{Bs}\,f(e^{-Bs}x,\mu)\,ds}\\[0.15in]
{\displaystyle
  =\lim_{T\rightarrow\infty}\,\frac{1}{T}\,\int_0^T\,e^{Bs}\,\left(D_xg(e^{-Bs}x,\mu)Be^{-Bs}x-Bg(e^{-Bs}x,\mu)\right)\,ds}\\[0.15in]
{\displaystyle
  =\lim_{T\rightarrow\infty}\,\frac{1}{T}\,\int_0^T\,\frac{d}{ds}\left(e^{Bs}g(e^{-Bs}x,\mu)\right)\,ds}\\[0.15in]
{\displaystyle
=\lim_{T\rightarrow\infty}\,\frac{e^{BT}g(e^{-BT}x,\mu)-g(x,\mu)}{T}}
\end{array}
\]
and this last limit is equal to $0$, since the numerator is bounded
in $T$ for any given $(x,\mu)\in\mathbb{R}^{\kappa+s}$.
So we conclude that $f\in\mbox{\rm ker}\,A$, and thus that
$\mbox{\rm ker}\,A=\mbox{\rm range}\,{\cal L}_{B}$.  This establishes (\ref{Aker}), and
concludes the proof of Proposition \ref{prop_Adecomposition}.
\hfill\qed

\end{document}